\documentclass[11pt]{article}
\usepackage[utf8]{inputenc}
\usepackage[T1]{fontenc}
\usepackage{amsmath,amssymb,amsfonts,amsthm}
\usepackage{physics} 
\usepackage{graphicx}
\usepackage{geometry}
\usepackage[colorlinks=true, pdfstartview=FitV, linkcolor=blue,citecolor=blue, urlcolor=blue]{hyperref}
\usepackage{xcolor}
\usepackage{authblk}
\usepackage[numbers,sort&compress]{natbib} 
\def\p{\partial}

\title{Resolving Sharp Gradients of Unstable Singularities \\
to Machine Precision via Neural Networks}

\author[1]{Yongji Wang}
\author[3]{Tristan L\'{e}ger}
\author[2,*]{Ching-Yao Lai}
\author[1,*]{Tristan Buckmaster}

\affil[1]{New York University, Department of Mathematics, New York, NY 10012, USA}
\affil[2]{Stanford University, Department of Geophysics, Stanford, CA 94305, USA}
\affil[3]{Yale University, Department of Mathematics, New Haven, CT 06511, USA}

\affil[*]{Corresponding authors: cyaolai@stanford.edu,  buckmaster@cims.nyu.edu}

\date{}

\begin{document}

\maketitle

\begin{abstract}
Recent work~\cite{WangEtAl2025Discovery} introduced a robust computational framework combining embedded mathematical structures, advanced optimization, and neural network architecture, leading to the discovery of multiple unstable self-similar solutions for key fluid dynamics equations, including the Incompressible Porous Media (IPM) and 2D Boussinesq systems. While this framework confirmed the existence of these singularities, an accuracy level approaching double-float machine precision -- i.e., accuracy constrained by inherent round-off errors machine precision ($O(10^{-13})$ residuals) -- was only achieved for stable and 1st unstable solutions of the 1D C\'ordoba-C\'ordoba-Fontelos model. For highly unstable solutions characterized by extreme gradients, the accuracy remained insufficient for validation. The primary obstacle is the presence of sharp solution gradients. Those gradients tend to induce large, localized PDE residuals during training, which not only hinder convergence, but also obscure the subtle signals near the origin required to identify the correct self-similar scaling parameter $\lambda$ of the solutions. In this work, we introduce a gradient-normalized PDE residual re-weighting scheme to resolve the high-gradient challenge while amplifying the critical residual signals at the origin for $\lambda$ identification. Coupled with the multi-stage neural network architecture, the PDE residuals are reduced to the level of round-off error across a wide spectrum of unstable self-similar singularities previously discovered. Furthermore, our method enables the discovery of new highly unstable singularities, i.e. the 4th unstable solution for IPM equations and a novel family of highly unstable solitons for the Nonlinear Schr\"{o}dinger equations. This results in achieving high-gradient solutions with high precision, providing an important ingredient for bridging the gap between numerical discovery and computer-assisted proofs for unstable phenomena in nonlinear PDEs.
\end{abstract}

\section{Introduction}

The study of singularity formation in fluids relies heavily on high-precision computation. Unstable singularities, especially in dimensions greater than 1, are exceptionally elusive; their discovery requires overcoming significant computational instabilities. Their rigorous mathematical validation via computer-assisted proofs (CAP) generally demands extremely high numerical accuracy \cite{Chen2021}. In this work, we refer to \emph{machine precision} as the highest achievable accuracy, where the remaining PDE residuals are dominated by the inherent round-off errors of double-float arithmetic, typically $O(10^{-13})$ or lower. While recent computational frameworks~\cite{WangEtAl2025Discovery} successfully identified families of unstable self-similar solutions in key fluid dynamics equations, achieving the requisite level of accuracy remains a profound computational challenge, particularly for highly unstable phenomena.

Highly unstable singularities are computationally pathological. As the order of instability increases, the solution profiles develop increasingly extreme localized gradients (see Figure 1). The framework introduced in Wang {\it et al.} \cite{WangEtAl2025Discovery}--built upon the earlier work \cite{Wang-Lai-GomezSerrano-Buckmaster:pinn-selfsimilar-boussinesq}--successfully discovered families of unstable solutions through the integration of mathematical analysis with customized machine learning techniques. Central to the methodology is an iterative feedback loop between numerical experiments and mathematical insights, which guide the design of specialized Physics-Informed Neural Network (PINN \cite{raissi2019physics,karniadakis2021physics}) architectures. Crucially, strong mathematical inductive biases are embedded directly into the network design. This includes input transformations to enforce symmetries, coordinate compactification to handle the infinite domains inherent in self-similar problems, and tailored \emph{solution envelopes} to ensure correct asymptotic decay and behavior near the origin. To achieve the extreme precision required, the framework moves beyond standard optimizers by employing a high-precision, full-matrix Gauss-Newton approach \cite{martens2015optimizing} and uses the Multi-Stage Neural Network (MSNN) \cite{Wang-Lai:multistage-nn} architecture. MSNN was successfully applied to the stable and 1st unstable solutions for CCF and IPM. Despite this, accuracy plateaued significantly for the most unstable examples. The precision of highly unstable solutions remained around $O(10^{-6})$ to $O(10^{-8})$, a level of precision insufficient for a MSNN method to effectively refine the solution.

The primary roadblock lies in the optimization landscape itself. In standard loss formulations, the  residuals generated in high-gradient regions overwhelm the optimization process. Crucially, this dominance obscures the subtle signals required to determine the exact physical parameters of the singularity, such as the self-similar scaling parameter $\lambda$ \cite{eggers2015singularities}. For the systems studied here, the Incompressible Porous Media (IPM) and Córdoba-Córdoba-Fontelos (CCF) equations, smoothness conditions at the origin (stagnation point) play a key role in determining the precise value of $\lambda$. The optimization stagnates, unable to reconcile the demands of the high-gradient regions with the precision required at the origin.

\begin{figure}[t]
	\centering
	\includegraphics[width=\linewidth]{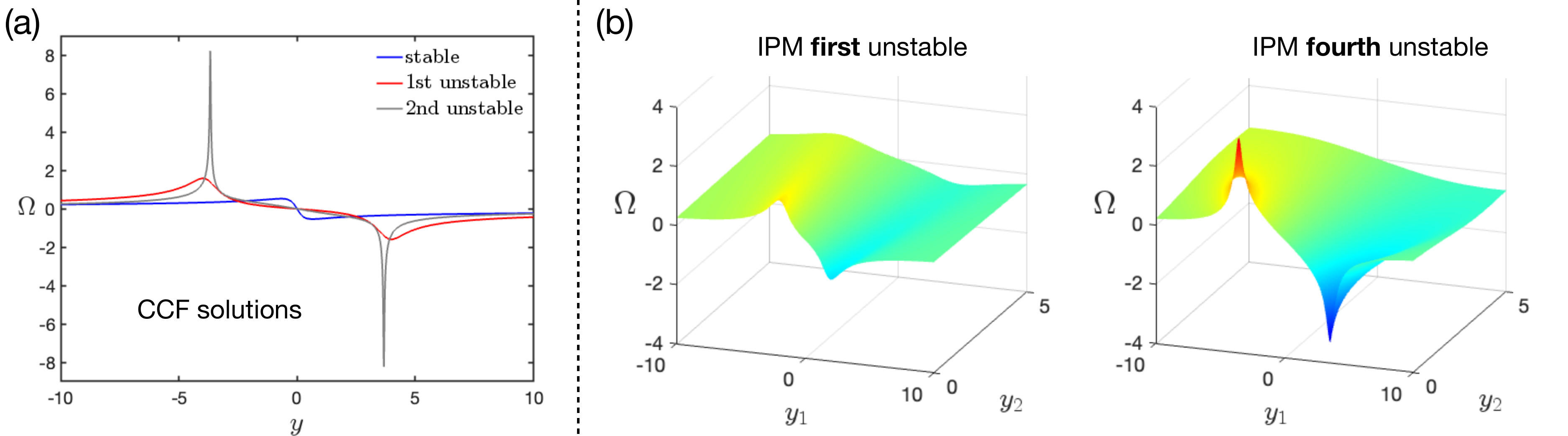}
	\caption{\small \textbf{High-gradient feature of unstable solutions} (a) Comparison of stable and unstable solution for (a) CCF equation and (b) IPM equations. The highly unstable solution shows significant higher local gradient than that of the stable or the low-order unstable solutions. }
	\label{fig:high-gradient}
\end{figure}

This paper introduces a critical innovation to overcome this barrier: a gradient-normalized loss function. By normalizing the PDE residuals relative to the local gradients of the solution, we effectively rebalance the optimization landscape. This approach prevents high-gradient regions from dominating the loss while simultaneously amplifying the critical residual signals at the origin, enabling the precise identification of $\lambda$. Additionally, integrating this gradient-normalization strategy with the established framework and MSNN training schemes, we demonstrate the ability to resolve these high-gradient, unstable singularities to near-machine precision. We achieve $O(10^{-13})$ residuals for the full suite of previously discovered CCF solutions and achieve accuracy between $O(10^{-11})$ and $O(10^{-13})$ for the 2-D IPM solutions. Furthermore, this enhanced precision enables the confirmation of the 4th unstable IPM solution.

The applicability of this high-precision framework extends beyond fluid singularities. We demonstrate its effectiveness in resolving unstable coherent structures in Nonlinear Schrödinger (NLS) equations, which share similar computational challenges, such as parameter sensitivity and infinite domains. We achieve double-float machine precision $O(10^{-13})$ in three such settings: unstable solitons in a double-well potential, Gross-Pitaevskii vortices with large winding numbers, and unstable excited states (solitons with nodes) of focusing NLS. The gradient renormalization also allows us to capture subtle local features of the solution that are not accessible through more standard PINN methods. For example, our approach reveals a novel empirical expression for the vortex core structure. For excited states, we identify a new empirical relationship between the azimuthal mode and the peak location.

\section{Machine-Precision Unstable Singularities}\label{sec:fluidSingularity}

The computational framework of \cite{WangEtAl2025Discovery} successfully identified families of unstable self-similar solutions in key fluid models. However, its accuracy plateaus for highly unstable solutions due to the emergence of extreme, localized gradients (Figure~\ref{fig:high-gradient}). These gradients create a pathological optimization landscape where large, localized PDE residuals overwhelm the subtle signals required to determine the self-similar scaling parameter, $\lambda$.

To overcome this fundamental barrier, we introduce a gradient-normalized residual. We demonstrate our method on two canonical models: the 1D C\'{o}rdoba-C\'{o}rdoba-Fontelos  equation and the 2D Incompressible Porous Media equation. By integrating this innovation with the architectural insights of \cite{WangEtAl2025Discovery}, our approach resolves these high-gradient solutions to machine precision $O(10^{-13})$
and confirms the fourth unstable solution for the IPM equation.

\subsection{Córdoba-Córdoba-Fontelos (CCF) Equation}\label{sec:ccf}
The CCF equation~\cite{Cordoba-Cordoba-Fontelos:CCF-model, Baker-Li-Morlet:analytic-structure-transport-equations} is a fundamental 1D model for nonlocal transport, serving as a critical analog for understanding singularity formation in more complex systems, such as the 3D Euler equations~\cite{Kiselev:open-problems}. We study the inviscid case:
\begin{align}\label{ccf}
\theta_t - H[\theta] \theta_x = 0,
\end{align}
where $H$ is the Hilbert transform. We seek self-similar blow-up solutions by analyzing the vorticity $\omega = \theta_x$, which is governed by the equation:
\begin{align}\label{ccf2}
 \omega_t-u\omega_x=\omega u_x,\qquad\mbox{where } u=\int_0^x(H\omega)(s)\,ds \, ,
\end{align}
where $u$ is considered the velocity of the system. To identify the self-similar profile of the blow-up solution as in \cite{WangEtAl2025Discovery,Wang-Lai-GomezSerrano-Buckmaster:pinn-selfsimilar-boussinesq}, we employ a self-similar ansatz
 \begin{equation}
 \omega=\frac{1}{1-t}\Omega\left(y\right)\quad \text{with} \quad y =\frac{x}{(1-t)^{1+\lambda}}, \quad \text{and} \quad U = \int_{0}^{y} H\Omega(s)\:\mathrm{d}s,\
 \end{equation}
for some $\lambda$ to be determined. Plugging that ansatz in \eqref{ccf2} yields the self-similar CCF equation:
  \begin{equation}\label{eq:ccf3}
 	 \Omega+((1+\lambda)y- U)\partial_y \Omega-\Omega \partial_y U=0\, .
 \end{equation}
Previous work successfully identified three self-similar solutions (stable \cite{Eggers_2020}, 1st unstable \cite{Eggers_2020,Wang-Lai-GomezSerrano-Buckmaster:pinn-selfsimilar-boussinesq}, and 2nd unstable \cite{WangEtAl2025Discovery}). While the first two solutions were resolved to machine precision \cite{WangEtAl2025Discovery}, the accuracy for the 2nd unstable solution was limited. As shown in Figure \ref{fig:high-gradient}(a), this solution exhibits significantly sharper gradients, posing a substantial computational challenge.

\begin{figure}[t]
	\centering
	\includegraphics[width=\linewidth]{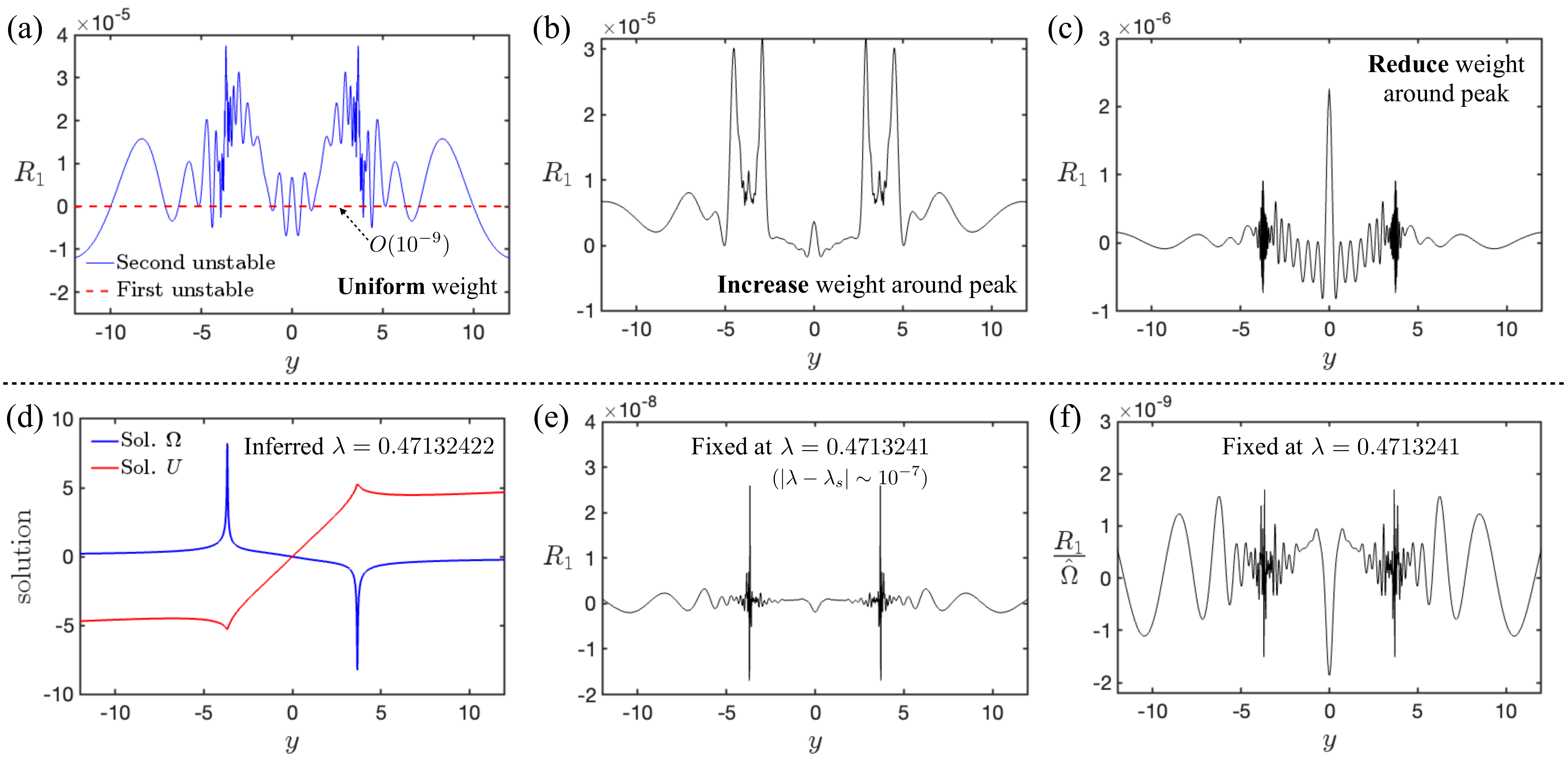}
	\caption{\small \textbf{Challenges of 2nd unstable solutions to CCF} (a) Comparison of the PDE residual for training the first and second unstable solutions for the CCF equation using the PINN framework with uniform weighting. (b \& c) The PDE residual of 2nd unstable solution trained by (b) increase or (c) reduce the weights of collocation points around the high-gradient area of the solution ($y=4$). (d) The spatial profile of the solutions $\Omega$ and $U$ for the 2nd unstable solution with the inferred $\lambda_s = 0.47132422$ by a single-stage training using gradient-normalized residual. (e \& f) The (e) absolute and (f) relative PDE residual for the CCF solution trained at the fixed $\lambda$ with $10^{-7}$ distance to the inferred smooth $\lambda_s$. } 
	\label{fig:ccf_issue}
\end{figure}

\subsubsection{Challenge of Learning High-Gradient Solutions}

The primary obstacle in resolving highly unstable solution for the CCF equation is the emergence of a pathological optimization landscape due to high gradients present in such solutions (Fig.~\ref{fig:high-gradient}a). In the absence of a technique to address this challenge, the training for the 2nd unstable solution stagnates, with the PDE residual plateauing at $O(10^{-5})$ (Fig.~\ref{fig:ccf_issue}a), orders of magnitude higher than the $O(10^{-9})$ achieved for the 1st unstable CCF solution (red dash line in Fig.~\ref{fig:ccf_issue}a). This high residual is dominated by contributions from the sharp gradient region near $y=4$. 

An intuitive remedy, such as increasing the number or weight of collocation points in this region, proves ineffective and can even degrade the overall solution quality (Figure \ref{fig:ccf_issue}b). Our systematic experiments revealed a counter-intuitive insight: {\it reducing} the loss contribution from the high-gradient region is necessary for convergence. As shown in Fig.~\ref{fig:ccf_issue}(c), manually down-weighting this region by approximately 10 times allows the residual to decrease by another order of magnitude before plateauing due to non-smoothness at the origin.

\subsubsection{Principle of Gradient-Normalized Residuals}

The underlying cause of this optimization failure is a mismatch between the training objective and the solution's properties. A standard PINN loss function, which minimizes the absolute PDE residual, implicitly drives the network to produce a solution with a uniform absolute error across the domain. For solutions with extreme gradients, this forces the optimization to excessively prioritize the high-gradient region while neglecting subtle but critical features elsewhere, such as at the origin. This imbalance ultimately causes the training to stagnate at a high equation residual, as shown in Figs.~\ref{fig:ccf_issue}(a) and (b).

The desired outcome, however, is a solution that exhibits a uniform relative error. To achieve this, we introduce our central contribution: a gradient-normalized residual. By normalizing the PDE residual in the loss function by the local solution gradient, we effectively re-balance the optimization landscape. This re-balancing ensures that the training process evaluates all regions of the domain equitably, regardless of the local gradient magnitude. As demonstrated in Fig.~\ref{fig:ccf_issue}(e), this approach enables the residual for the 2nd unstable solution to converge to a level comparable to that of the first unstable solution. Furthermore, the gradient-normalized residual amplifies the signal at the origin, preventing it from being obscured by the large residual in the peak region (Fig.~\ref{fig:ccf_issue}f). This capability is critical for the precise identification of the self-similar scaling parameter $\lambda$. The detailed justification and implementation of the gradient-normalized residual approach are given in Section \ref{sec:grad_method}.

\begin{figure}[t]
	\centering
	\includegraphics[width=\linewidth]{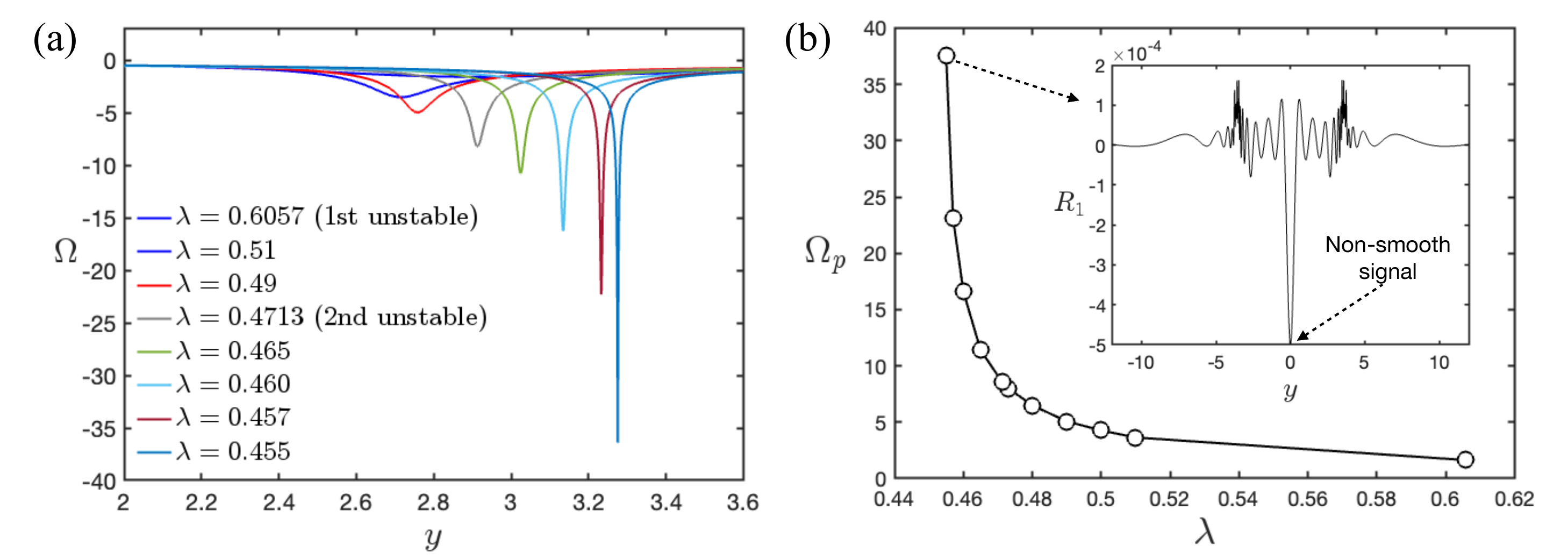}
	\caption{\small \textbf{CCF solutions with low $\lambda$.} (a) Comparison of smooth and non-smooth CCF solutions for $\lambda < 0.61$. The magnitude and the sharpness of the high-gradient areas increase with the decrease of $\lambda$ (b). Evolution of the peak value of the CCF solution with respect to $\lambda$, illustrating the intensification of gradients as instability increases. The inset shows that PDE residual for the CCF solution at the fixed $\lambda = 0.455$, which has non-smooth signal at the origin, indicating that it remains a non-smooth solution.}
	\label{fig:CCF_high_peak}
\end{figure}

\subsubsection{Exploring higher Unstable Solutions}
Having confirmed the effectiveness of the gradient-normalized residual approach, we extended our search to the third unstable solution. Our investigation shows that as $\lambda$ decreases, the solution's peak amplitude increases dramatically, reaching nearly 40 for $\lambda\approx 0.455$ (Fig.~\ref{fig:CCF_high_peak}). This suggests that a potential third unstable solution would possess even more extreme gradients than the second.

While our gradient-normalization method effectively reduces the residual in the peak region, a persistent non-smooth signal at the origin (inset of Fig.\ref{fig:CCF_high_peak}b) remains for all tested $\lambda \in[0.455, 0.4713]$. This indicates that if a third unstable solution exists, it is located at a $\lambda$ value corresponding to an even higher-gradient profile. The challenge of exploring the range $\lambda < 0.455$ and resolving such a feature in the 1D CCF model appears to be even greater than finding highly unstable solutions in the 2D IPM equation. A full investigation into the existence and properties of this third unstable solution is a compelling direction for future research.

\subsection{Incompressible Porous Media (IPM) Equation}\label{sec:ipm}

\begin{figure}[t!]
	\centering
	\includegraphics[width=\linewidth]{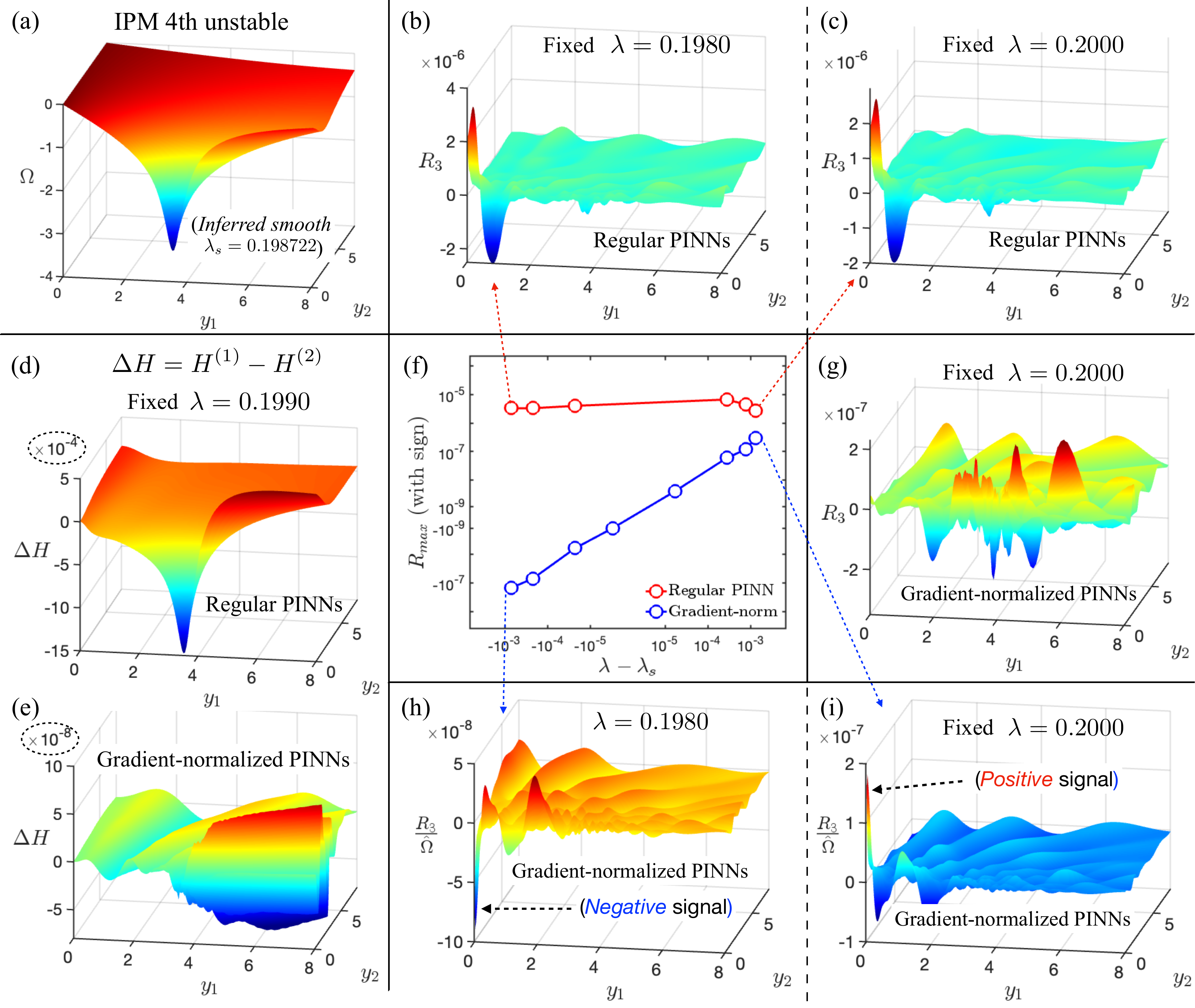}
	\caption{\small \textbf{Challenge of 4th unstable solution to IPM.} (a) The vorticity profile $\Omega$ for the 4th unstable solution to IPM with the inferred smooth $\lambda_s$. (b~\&~c) The PDE residual from the PINN training using standard equation loss at two different $\lambda$ around the 4th unstable solution, where the spurious non-smooth signals at the origin are nearly identical. (d) The difference between the solutions $H^{(1)}$ and $H^{(2)}$, which are trained at the same fixed $\lambda = 0.1990$ via the same PINN framework using standard equation loss but different random initialization of the network weight. The solution discrepancy $O(10^{-4})$ is much larger than the PDE residual $O(10^{-6})$ (panel b) for each solution. (e) In comparison, the difference between the solution using gradient-normalized residual remains the same order of magnitude with the PDE residual $O(10^{-8})$ (panel h). (f) The relation of the maximum value $R_{max}$ of the relative PDE residual at the origin captured via PINN framework using standard equation loss (red) and gradient-normalized residual (blue) with the distance of $\lambda$ to the true smooth $\lambda_s$. (g) the absolute PDE residual obtained via PINN framework using gradient-normalized residual. (h~\&~i) The relative PDE residual normalized by the solution magnitude at two different $\lambda$, which reveals the true non-smooth signal at the origin that are proportional to the distance from the true smooth $\lambda$ (panel e).}
	\label{fig:IPM4th}
\end{figure}

The 2D IPM equation models the dynamics of an incompressible fluid in a porous medium governed by Darcy's law~\cite{Bear:dynamics-porous-media}. We analyze the system on the upper half-plane with non-penetration boundary conditions:
\begin{equation}\label{eq:IPM_time}
    \p_t \rho+ {\bf u}\cdot \nabla \rho=0,\quad \nabla \cdot {\bf u}=0,\quad {\bf u}+\nabla p=-(0,\rho),
\end{equation}
where $p$ is the pressure and ${\bf u}=(u_1, u_2)$ is the velocity vector with $u_2$ satisfying the no-penetration condition, namely $u_2(x_1, x_2=0) = 0$. Analogous to the CCF equation, we identify the self-similar blowup solution of IPM by utilizing the self-similar ansatz
 \begin{equation}\label{eq:IPM_ansatz}
 \rho({\bf x}, t)=(1-t)^\lambda H\left({\bf y}\right)\quad \text{and} \quad {\bf u}({\bf x}, t) = (1-t)^\lambda {\bf U}\left({\bf y}\right) \qquad \text{with} \quad {\bf y} =\frac{\bf x}{(1-t)^{1+\lambda}}.
 \end{equation}
where ${\bf x}=(x_1, x_2)$ are the original 2-D Cartesian coordinates, and ${\bf y}=(y_1, y_2)$ are the self-similar coordinates. Taking the curve of the third equation in \eqref{eq:IPM_time}, we can eliminate the pressure $p$ in the equation. Then, substituting the ansatz \eqref{eq:IPM_ansatz} gives the self-similar equations of IPM: 
\begin{gather}\label{IPM_self}
 -\lambda H + [(1+\lambda) y + {\bf U}] \cdot \nabla H = 0 \\
 \partial_{y_1} U_1 + \partial_{y_2}U_2 = 0  \qquad \text{and} \qquad 
 \partial_{y_1} U_2 - \partial_{y_2}U_1 = \partial_{y_1} H
\end{gather}
Wang {\it et al.} \cite{WangEtAl2025Discovery} identified, for the first time, four self-similar solutions, ranging from the stable solution to the third unstable solution, via PINNs with a particular design of coordinate compactification and a structured solution ansatz to enforce the symmetry and asymptotic behavior of the solutions (see Section \ref{sec:basicFW}). However, this new framework shows limitations in finding the fourth unstable solution, which also exhibits a sharp gradient (Fig.~\ref{fig:IPM4th}a), similar to the 2nd unstable CCF solution.

\subsubsection{Spurious Non-Smooth Signal induced by High Gradients} 
For low $\lambda$, the high-gradient nature of the IPM solution poses the same computational challenges as those encountered in the CCF problem. When using the PINN frameworks with uniform weighting of collocation points across the domain, the PDE residual for various fixed $\lambda$ values near the 4th unstable solution consistently shows a maximum residual at the origin, seemingly capturing the correct non-smooth signal of the solution (Fig.~\ref{fig:IPM4th}b). However, a closer comparison (Fig.~\ref{fig:IPM4th}b~\&~c) reveals that these signals at the origin are nearly identical across the range of $\lambda\in[0.1980, 0.2]$ (red line in Fig.~\ref{fig:IPM4th}e). This lack of variation prevents the identification of the correct signal needed to determine the smooth $\lambda$ for the fourth unstable solution. 

This erroneous signal at the origin is a direct consequence of the high gradient feature of the solution in the nearby region. We find that, despite the high-gradient profile of the solution $\Omega$ in the region around $(y_1, y_2) \approx (3.5, 0)$ (Fig.~\ref{fig:IPM4th}a), the magnitude of the PDE residual there is comparable to that in the rest of the domain (Fig.~\ref{fig:IPM4th}b). This indicates that the optimization process has excessively prioritized minimizing the PDE residual in that specific region, which could induce a risk of eliminating important solution features in other regions. To confirm this, we conducted multiple training experiments at a fixed $\lambda=0.1990$ using identical PINN settings but different random initializations of network weights. While all experiments converged to a PDE residual of $O(10^{-6})$ (similar to Fig.~\ref{fig:IPM4th}b), the average difference between any two trained solutions reaches $O(10^{-4})$, with discrepancies in the high-gradient region reaching $O(10^{-3})$ (Fig.~\ref{fig:IPM4th}d). This demonstrates that while the training improves the solution locally in the high-gradient area, it fails to capture the global features consistently, which ultimately prevents the optimization from converging to the expected unique solution. Consequently, the large-residual signal at the origin in Fig.~\ref{fig:IPM4th}(b) fails to reflect the true non-smooth behavior of the solution.

\subsubsection{Application of Gradient-normalized residual}
As with the CCF equation, the remedy for the high-gradient challenge is to incorporate the gradient-normalized residual into the loss function, ensuring that the training process can evaluate all regions of the domain equitably. Figure \ref{fig:IPM4th}(g) shows the PDE residual obtained at the same fixed $\lambda = 0.20$ as in Fig.~\ref{fig:IPM4th}(c), but it is now obtained using the gradient-normalized approach. The training achieves a significantly lower error across the entire domain, including the high-gradient regions. The incorrect non-smooth signal previously seen in Fig.~\ref{fig:IPM4th}(b) is no longer present (Fig.~\ref{fig:IPM4th}g); instead, the true non-smooth signal is revealed by plotting the relative PDE residual normalized by the solution magnitude (Fig.~\ref{fig:IPM4th}h and i), which isolates the true signal from the otherwise large residuals in the high-gradient areas.

To confirm that the gradient-normalized residual approach yields a unique solution, we repeated the experiments at the same fixed $\lambda=0.1990$ with different random initializations of network weights, applying the gradient-normalized residual approach this time. We clearly see that the difference between any two experiments is now on the order of $10^{-8}$ (Fig.~\ref{fig:IPM4th}e), which is consistent with the magnitude of the PDE residual (Fig.~\ref{fig:IPM4th}h). This confirms that the training now converges to a unique solution profile.

Finally, to validate the non-smooth signal identified by our new gradient-normalized approach, we recorded the magnitude and sign of the normalized PDE residual $\hat{R_3}=R_3/\hat{\Omega}$ at the origin (Fig.~\ref{fig:IPM4th}h~\&~i) for $\lambda$ values at varying distances from the true smooth value $\lambda_s$. Crucially, we found that the sign of this signal at the origin flips as the fixed $\lambda$ crosses the true smooth value $\lambda_s$ (Fig.~\ref{fig:IPM4th}h~\&~i). Furthermore, the magnitude of the signal is proportional to the distance from the true smooth $\lambda_s$, namely $|\lambda-\lambda_s|\propto \hat{R}_3(0)$, as shown in Fig.~\ref{fig:IPM4th}(e). This relationship successfully recovers the characteristic funnel plot previously observed for the low-order unstable solutions in Wang {\it et al.} (2025)\cite{WangEtAl2025Discovery}, confirming the robustness of our new method.

\subsection{Achieving Machine Precision}
The gradient-normalized residual allows us to resolve highly unstable singularities to an accuracy comparable to that of less unstable modes. To further increase this accuracy, we must also employ multistage training. This section describes how we integrate the gradient-normalized residual approach into the multistage training scheme, establishing a robust methodology to resolve unstable singularities to machine precision. Furthermore, we demonstrate how this combined approach enables high-precision identification of the self-similar scaling parameter $\lambda.$

\subsubsection{High-gradient challenge on multistage training}\label{sec:multiGrad}
\begin{figure}[!t]
	\centering
	\includegraphics[width=\linewidth]{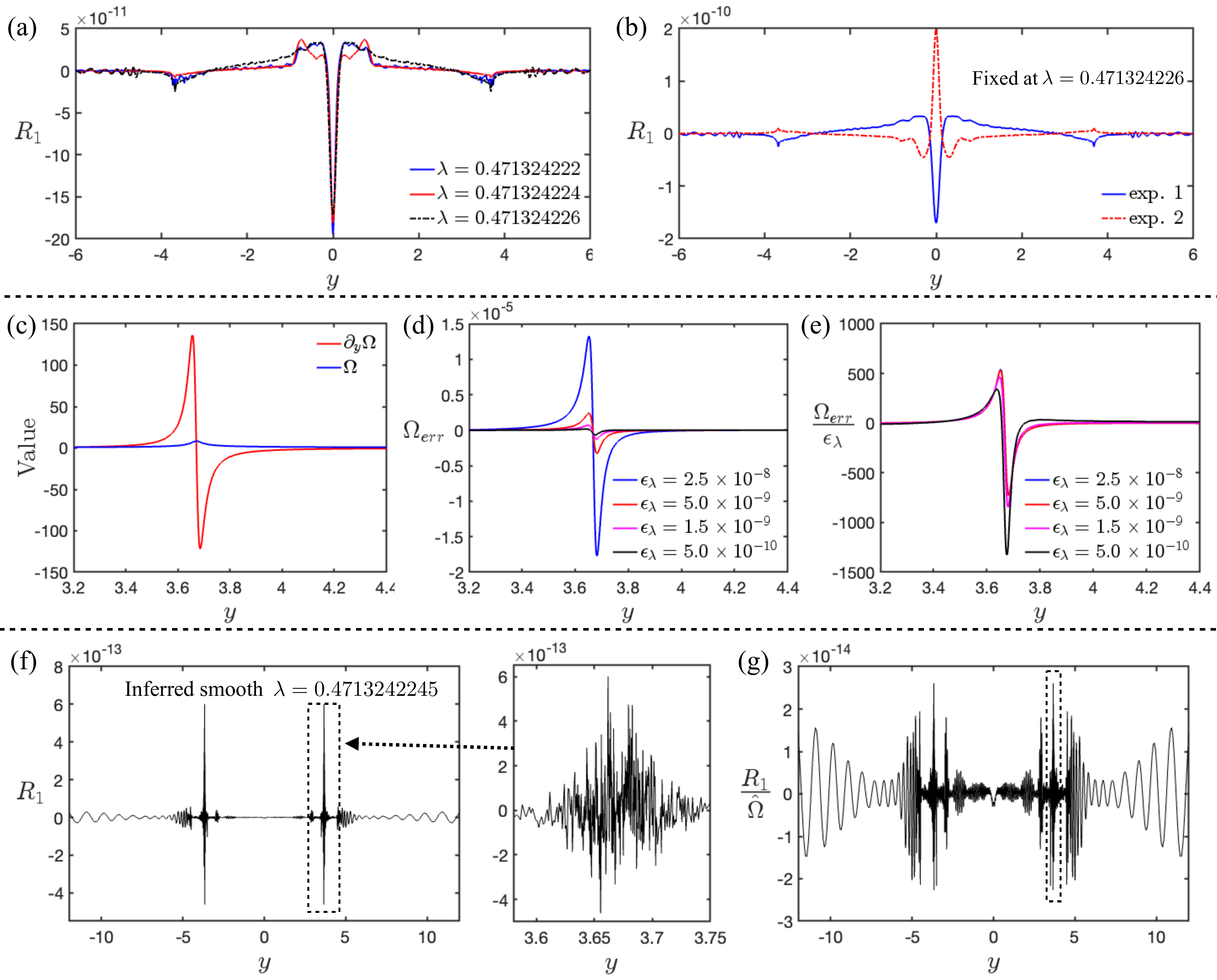}
	\caption{\small \textbf{Multistage training on 2nd unstable CCF solutions} (a) PDE residuals after second-stage training using the gradient-normalized residual based on the solution magnitude, which provide spurious non-smooth signal that remains identical for different $\lambda$ (b) Spurious signals for a fixed $\lambda$ are inconsistent, flipping sign between experiments with different random initializations (c) Comparison between the profiles of $\Omega$ and its derivative $\partial_y\Omega$. The latter has much higher gradient. (d) The solution error $\Omega_{err} = \Omega_g - \Omega_0$ between the {\it ground-truth} solution and the first-stage results $\Omega_0$ trained at different fixed $\lambda$. $\epsilon_\lambda$ indicates the proximity of the fixed $\lambda$ to the {\it ground-truth} $\lambda_s$ for the 2nd CCF unstable solution. Here, the ground-truth solution $\Omega_g$ and lambda $\lambda_s$ refer to the best 2nd CCF unstable solution obtained that reaches machine precision as shown in e.  (e) The normalized solution error by the error $\epsilon_\lambda$ of $\lambda$. (f) The absolute PDE residual for the 2nd CCF unstable solution using the new gradient-normalized residual based on the derivative of the solution $\partial_y\Omega_0$. The peak in the region of high solution gradient is restricted by the round-off error of double-float precision. (g) The relative PDE residual normalized by the solution magnitude.}
	\label{fig:CCFmulti}
\end{figure}

After obtaining an accurate solution from the first stage, a new challenge emerges. Figure \ref{fig:CCFmulti}(a) shows the PDE residuals for the second unstable CCF solution after applying the multistage method, which also uses the gradient-normalized residual. Compared with the first-stage results (Fig.~\ref{fig:ccf_issue}e), the residuals across most of the domain, including the high-gradient region (around $y=4$), are reduced by two to three orders of magnitude. However, the training stagnates, revealing a maximum PDE residual at the origin (Fig.~\ref{fig:CCFmulti}a) that is far larger than the residuals in the rest of the domain.
This signal seemingly captures the solution's non-smooth feature at the origin. However, it is spurious. By comparing the PDE residual for different $\lambda$ values in the narrow range $\lambda\in[0.471324222,0.471324226]$, there is no meaningful variation when comparing PDE residuals (Fig.~\ref{fig:CCFmulti}a). 
This resembles the spurious signal encountered when solving for the 4th IPM unstable solution (Fig.~\ref{fig:IPM4th}b) and prevents us from finding a more accurate smooth $\lambda_s$. To further confirm this, we conducted experiments at the exact same $\lambda$ but with different random initializations for the second-stage network. We found that the sign of the maximum residual flips between experiments, demonstrating it is a numerical artifact rather than a reliable indicator of the true non-smooth signal. 

Here, we note a critical distinction: the optimization artifact in the IPM case (Fig.~\ref{fig:IPM4th}b) arose from a standard equation loss; whereas, in the CCF second-stage training (Fig.~\ref{fig:CCFmulti}a), a similar artifact persists despite the use of the gradient-normalized residual. This demonstrates that the first-stage residual normalization is insufficient for the second-stage training. The reason is that the target of the second stage, the error function $\Omega_e$, possesses an even sharper gradient than the solution $\Omega_0$ itself. To see this, we introduce the solution ansatz for the multistage training:
\begin{equation}
    \Omega_g = \Omega_0 + \epsilon \Omega_e, \quad  U_g = U_0 + \epsilon U_e, \quad  \text{and} \quad \lambda_s^{(2)} = \lambda_0 + \epsilon_\lambda
\end{equation}
where $\Omega_0$, $U_0$, and $\lambda_0$ represent the solutions and $\lambda$ learned in the first-stage training. $\Omega_g$ and $U_g$ are the exact solutions of the 2nd CCF unstable solution at the exact smooth $\lambda_s^{(2)}$. Here, $\epsilon$ and $\epsilon_\lambda$ denote the magnitude of the error in the first-stage solutions $(U_0, \Omega_0)$ and $\lambda_0$. $\Omega_e$ and $U_e$ are the normalized error functions of the first-stage solutions, which are the target functions for the second-stage training. Assuming the first stage is accurate, we substitute the ansatz into the self-similar CCF equation \eqref{eq:ccf3} and neglect higher-order terms ($o(\epsilon)$ and $o(\epsilon_\lambda)$), which yields a linearized equation for the second-stage training
\begin{equation}\label{eq:ccf2ndeqn}
    \Omega_e+((1+\lambda)y- U_0)\partial_y \Omega_e - U_e\partial_y \Omega_0 -\Omega_e \partial_y U_0 - \Omega_0 \partial_y U_e=-\frac{R_1^{(0)}}{\epsilon} - \frac{\epsilon_\lambda}{\epsilon}\partial_y\Omega_0\, 
\end{equation}
The two terms on the right-hand side are the source terms, which are known from the first stage. The first term arises from the first-stage PDE residual $R_1^{(0)}$, while the second term arises from the error in the first-stage $\lambda$. 

We note that the second source term depends solely on the derivative of $\Omega_0$, which has a much larger magnitude and sharpness than $\Omega_0$ itself (Fig.~\ref{fig:CCFmulti}c). Thus, when the error from $\lambda_0$ is dominant (i.e., $\epsilon_\lambda$ is large relative to the residual $R_1^{(0)}$), the error function $\Omega_e$ will be dominated by the profile of $\partial_y\Omega_0$, not $\Omega_0$. Figure \ref{fig:CCFmulti}(d) shows the solution error $\Omega_{\text{err}} = \Omega_g - \Omega_0$ between the ground-truth profile and the first-stage results. When the error in $\lambda_0$ ($\epsilon_\lambda = \lambda_s-\lambda_0$) is large, the magnitude of $\Omega_{\text{err}}$ is proportional to $\epsilon_\lambda$, and its profile indeed approximates that of $\partial_y\Omega_0$. After normalizing the solution error by $\epsilon_\lambda$, all the profiles collapse onto the same curve (Fig.~\ref{fig:CCFmulti}d).

Since the solution profile for the second-stage training, $\Omega_e$, is proportional to $\partial_y\Omega_0$, we must normalize the equation residual by a new function that reflects the sharp gradient of $\partial_y\Omega_0$, rather than $\Omega_0$. In addition, as the magnitude of $\Omega_e$ is extremely large in the high-gradient region (Fig.~\ref{fig:CCFmulti}e), we multiply the second-stage network output by an \textit{envelope} that captures this magnitude. By integrating these two approaches into the multistage training, we successfully achieved machine precision ($O(10^{-13})$ residuals) for the 2nd CCF unstable solution (Fig.~\ref{fig:CCFmulti}f). We note that although the absolute PDE residual remains higher in the region of high solution gradient, this is constrained by the round-off error of double-float precision, as expected for a large solution magnitude there. The relative PDE residual, normalized by the solution magnitude, achieves uniform machine precision (Fig.~\ref{fig:CCFmulti}g).

This enhanced multistage training framework can also be applied to the 2D IPM problem, which significantly improves the accuracy of all unstable solutions. More detailed results will be discussed in later sections.

\subsubsection{High-precision identification of $\lambda$}
Discovering an unstable singularity requires not only finding the solution with high accuracy but also identifying the correct value of the scaling parameter $\lambda$ corresponding to a smooth solution. Therefore, to confirm the effectiveness of the enhanced multistage method as described above, we must verify that it can precisely identify $\lambda$, in addition to reducing the PDE residuals to round-off error (Fig.~\ref{fig:CCFmulti}g). Here, we note that all results described here utilize the enhanced PINN framework with gradient-normalized residuals. The {\it non-smooth signal} discussed here refers to the maximum value $R_{max}$ of the relative PDE residual at the origin, defined as the absolute PDE residual normalized by the solution magnitude (Fig.~\ref{fig:IPM4th}h).

\begin{figure}[t]
	\centering
	\includegraphics[width=\linewidth]{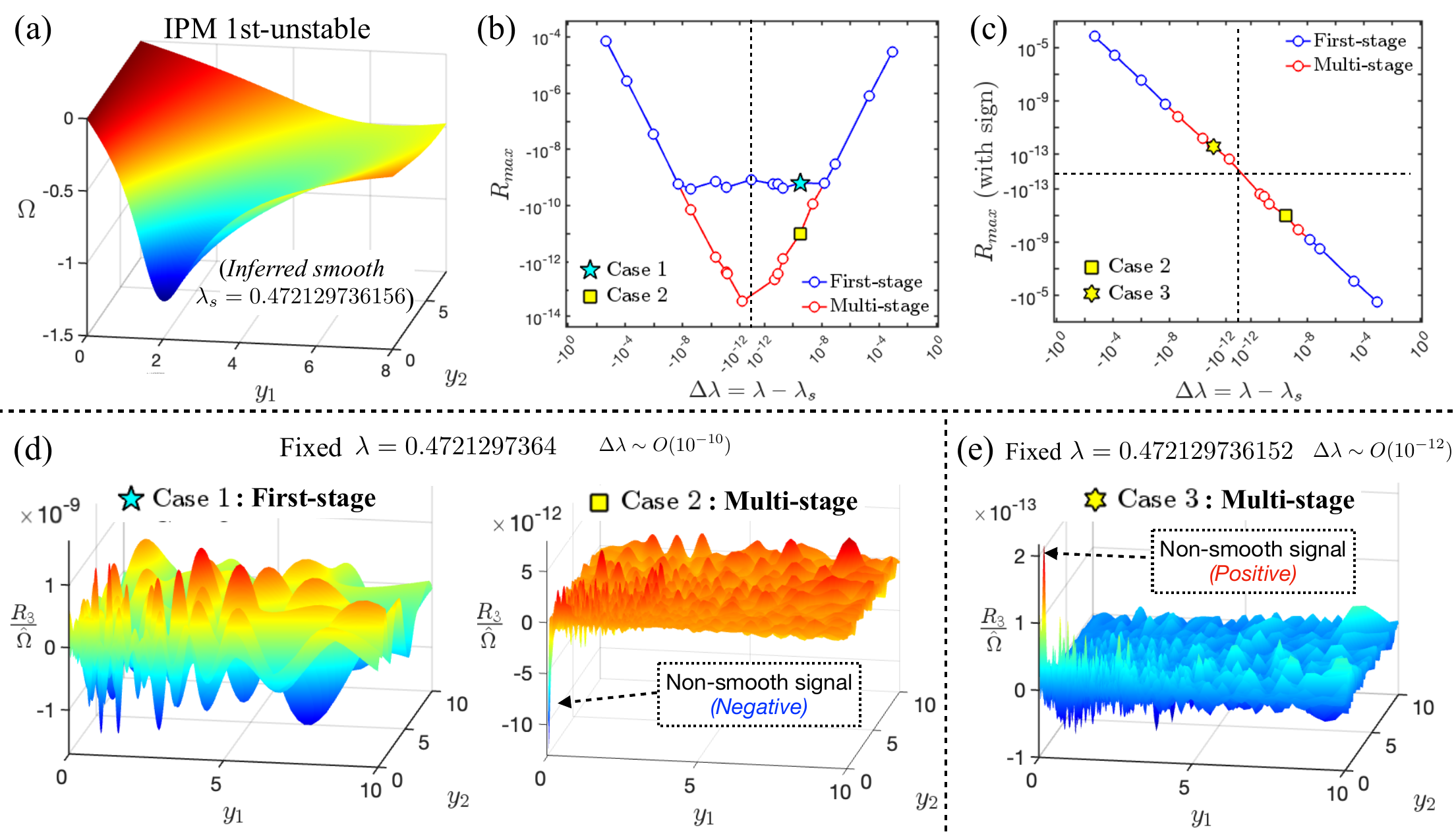}
	\caption{\small \textbf{Multistage training on IPM unstable solution}. (a) The vorticity profile $\Omega$ for the 1st unstable solution to IPM equation with the inferred smooth $\lambda_s$. (b) The relationship between the magnitude of maximum signal $R_{max}$ at the origin and the distance $\Delta \lambda = \lambda-\lambda_s$ to the smooth $\lambda_s$, obtained via single-stage training (blue line) and multistage training (red line). (c) The relationship between the signed $R_{max}$ and the distance $\Delta \lambda$ follows a perfect straight line on a log-log space. (d) Comparison of the relative PDE residual obtained via single-stage training and multistage training at the same fixed $\lambda$ with $\Delta \lambda \sim O(10^{-10})$. The latter one reveals the true non-smooth signal at the origin with a negative sign. (e) The relative PDE residual obtained via multistage training at the fixed $\lambda$ with $\Delta \lambda \sim O(10^{-12})$, which reveals the non-smooth signal with a positive sign.}
	\label{fig:2DlambMulti}
\end{figure}

The blue line in Fig.~\ref{fig:2DlambMulti}(b) illustrates the relationship between the magnitude of the non-smooth signal at the origin and the distance $\Delta=\lambda-\lambda_s$, obtained via single-stage training. Using the 1st IPM unstable solution as an example (Fig.~\ref{fig:2DlambMulti}a), for $|\Delta \lambda| > 10^{-8}$, we observe a clear power-law relationship between $R_{max}$ and $\Delta\lambda$, forming the characteristic {\it funnel plot} as reported in Wang {\it et al.}(2025). However, the accuracy of single-stage training is limited; the average relative PDE residual stagnates around $10^{-9}$ (Fig.~\ref{fig:2DlambMulti}d). Consequently, when $|\Delta\lambda|<10^{-8}$, the expected true signal $R_{max}$ falls below this background residual level, which masks the signal and prevents further refinement of $\lambda$. This also causes the observed $R_{max}$ to plateau at the average PDE residual level $O(10^{-9})$ (blue line in Fig.~\ref{fig:2DlambMulti}b).

Therefore, minimizing the overall PDE residuals is essential for resolving the non-smooth signal required for high-accuracy identification of $\lambda$. Applying the enhanced multistage method (Section \ref{sec:multiGrad}), we conduct the second-stage training on the IPM solution at a fixed $\lambda$ with $|\Delta\lambda| \sim 10^{-10}$. The multistage training successfully reduced the average PDE residual across the domain by three orders of magnitude (Fig.~\ref{fig:2DlambMulti}d). Crucially, the true non-smooth signal at the origin re-emerged as the maximum residual in the domain, continuing its expected decay $O(10^{-11})$ based on the given $|\Delta\lambda|$. Figure \ref{fig:2DlambMulti}(e) further shows results for $|\Delta\lambda| \sim 10^{-12}$; the multi-stage training still resolves the background residual and reveals the true non-smooth signal at the origin, with the magnitude reduced to $O(10^{-13})$.

Moreover, systematic experiments using the multistage method across a range of $\lambda$ successfully extended the funnel plot relationship between $R_{max}$ and $\Delta \lambda$ down to $O(10^{-13})$, close to double-float machine precision (Fig.~\ref{fig:2DlambMulti}b, red line). The extended relationship follows the same power law identified in the first stage. Furthermore, the sign of the non-smooth signal $R_{max}$ correctly indicates whether the fixed $\lambda$ is larger or smaller than the true smooth $\lambda_s$ (Fig.~\ref{fig:2DlambMulti}d vs e). Plotting the signed $R_{max}$ against $\delta_\lambda$ reveals a perfect straight line on a log-log plot (Fig.~\ref{fig:2DlambMulti}c). This further confirms the authenticity of the refined non-smooth signal resolved by the multistage method and establishes a quantitative link between the maximum residual at the origin and the accuracy of the inferred $\lambda$.

\subsubsection{Comparison of $\lambda$ signals across IPM singularities}
\begin{figure}[t!]
	\centering
	\includegraphics[width=\linewidth]{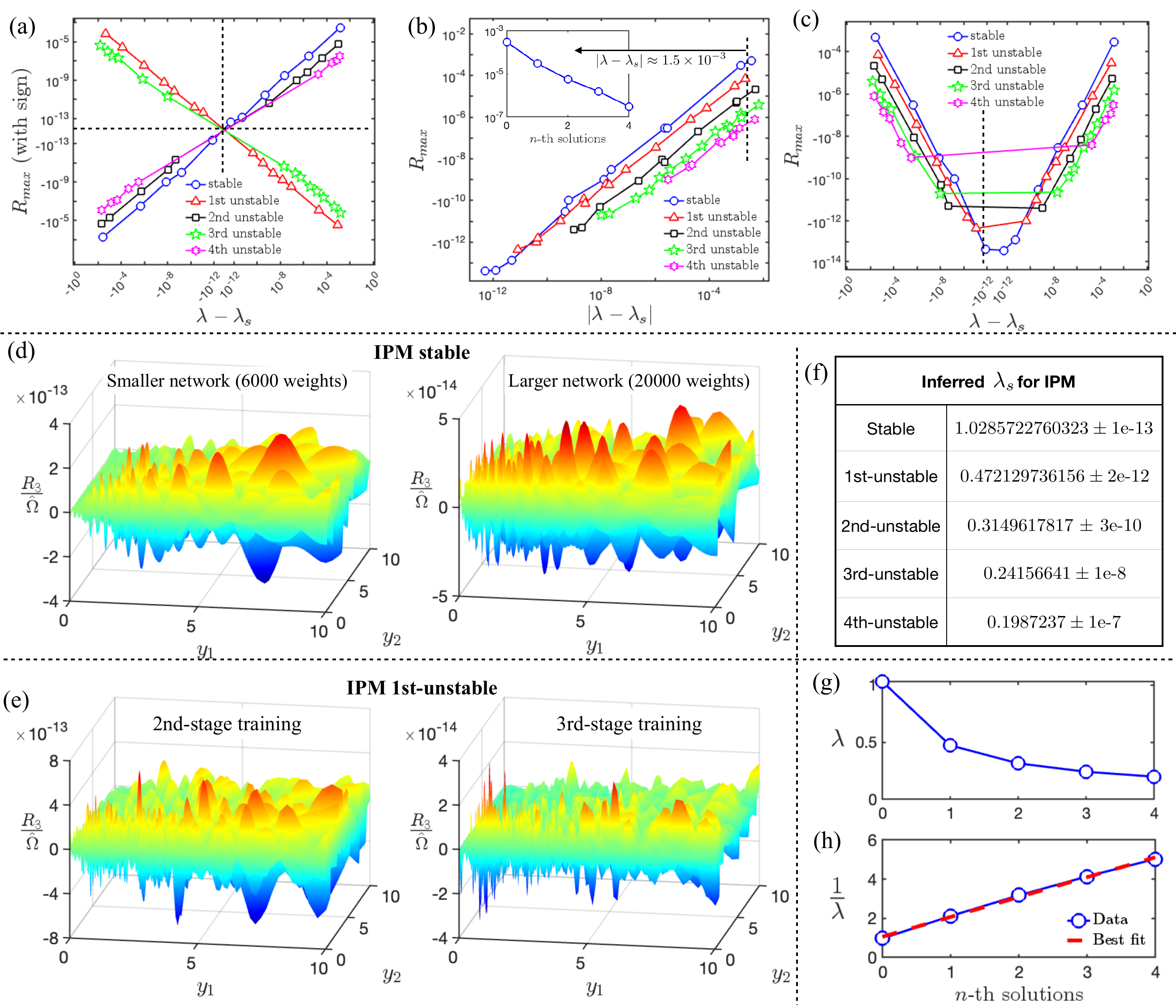}
	\caption{\small \textbf{High-precision $\lambda$ identification} (a) The relationship between the signed non-smooth signal $R_{max}$ and the signed distance $\Delta\lambda=\lambda-\lambda_s$ for different IPM singularities. (b) The relationship between the magnitude of $R_{max}$ and the absolute distance $|\Delta\lambda|$. The inset shows that for the same distance $|\Delta\lambda|$, the non-smooth signal $R_{max}$ detected at the origin decays exponentially for higher-order modes. (c) The horizontal line implies the lowest non-smooth signal $R_{max}$ for each mode that can be detected using a fixed computational budget with two stage of training. (d) Comparison of the relative PDE residual for the stable IPM singularity via two-stage training using networks of different sizes (i.e. number of weights). (e) Comparison of the relative PDE residual for the 1st IPM unstable solution using two stages and three stages of training. (f) The list of the inferred smooth $\lambda_s$ and their accuracy for different IPM singularities. (g) The empirical relation of the smooth $\lambda_s$ with the order $n$ of the mode. (h) The relation of $1/(\lambda-\lambda_0)$ versus $n$ falls on a nearly perfect straight line, with $\lambda_0=0.031694$ obtained by the best fit. }
	\label{fig:2DlambSignal}
\end{figure}

Our enhanced multistage method successfully resolved all IPM singularities, from the stable mode $(n=0)$ to the 4th unstable mode $(n=4)$, to high precision. This allows us to compare the properties of the non-smooth signal $R_{max}$, used to identify $\lambda_s$, among different solutions. Figure \ref{fig:2DlambSignal} shows the relationship between the signed signal $R_{max}$ and the distance $\Delta \lambda$ for all singularities. We find that all of them follow a straight line in the log-log space, but the sign of the slope alternates for each successive mode. Namely, for the stable solution, a positive non-smooth signal indicates that the current $\lambda$ is higher than the smooth one; however, for the 1st unstable solution, a positive signal indicates a lower $\lambda$.

Besides the alternating slope sign, the curve relating the magnitude of $R_{max}$ and the absolute distance $|\Delta_\lambda|$ is also shifted downward for higher-order singularities (Fig.~\ref{fig:2DlambSignal}b). For example, with the $\lambda$ distance around $O(10^{-3})$, the non-smooth signal detected at the origin for the stable solution is around $O(10^{-3})$, while for the 4th unstable solution, the non-smooth signal becomes $O(10^{-7})$. On average, the magnitude of the non-smooth signal $R_{max}$ decays exponentially, reducing by one order of magnitude for each higher-order singularity. This feature confirms that resolving higher-order modes is intrinsically more difficult, as an exponentially weaker $R_{max}$ signal demands a progressively lower background PDE residual to be revealed. In this case, the multistage method is essential for finding higher-order unstable singularities with precision.

While our enhanced multi-stage method is robust, the high-gradient nature of higher-order unstable modes still reduces the convergence rate of training. For a fixed computational budget (identical network size, iterations, and two training stages), the achievable background PDE residual for the smooth solution increases by approximately one order of magnitude for each successive singularity, which limits the non-smooth signal $R_{max}$ to be detected. For example,  figure~\ref{fig:2DlambSignal}(c) shows that the best $R_{max}$ residual for the stable solution reaches $O(10^{-14})$, while the one for the 4th unstable solution stagnates at $O(10^{-10})$. This increased background residual, combined with the exponentially weaker $R_{max}$ signal, decreases the achievable $\lambda$ accuracy for a fixed computational budget by approximately two orders of magnitude per model. 

However, the effectiveness of multistage training can be increased by using larger networks with more weights (Fig.~\ref{fig:2DlambSignal}d), or by introducing a third stage to further reduce the PDE residual and refine the $\lambda$ identification (Fig.~\ref{fig:2DlambSignal}d). Both approaches, however, require more memory. Figure \ref{fig:2DlambMulti}(f) shows the best inferred $\lambda$ for different singularities of IPM using a multistage setup that can be implemented on a single A100 GPU (80GB). To further improve the $\lambda$ accuracy, one could use larger networks or more training stages in parallel on multiple GPUs.

With high-precision $\lambda$ values for five IPM singularities, we re-examine the empirical relation between the smooth $\lambda_s$ and the order $n$ of the mode (Fig.~\ref{fig:2DlambSignal}g). Here, we adopt the same relation ansatz as used in Wang {\it et al.} (2025) \cite{WangEtAl2025Discovery}
\begin{equation}\label{eq:lamblaw}
    \lambda = \frac{1}{a + b\lambda}
\end{equation}
with $a$, $b$ as unknown parameters to be determined. The discovery of the 4th unstable solution provides an additional data point of $\lambda_s$ to check the robustness of this relation. By linear regression, we found the best-fit coefficients in \eqref{eq:lamblaw} to be $a=1.014$, $b=1.059$ with a coefficient of determination ($r^2$ value) reaching 0.997, indicating a near-perfect fit. Figure \ref{fig:2DlambSignal}(h) shows that the plot of $1/\lambda$ versus $n$ falls on a perfectly straight line.

\section{Application to Unstable Coherent Structures}
The high-precision framework developed in this work extends beyond the fluid dynamics equations studied above and is applicable to a wide range of problems involving unstable solutions and coherent structures defined on infinite domains. Nonlinear Schr\"{o}dinger equations (NLS) are prominent examples that present similar numerical challenges, requiring precise computation to capture elusive solutions. In this study, we focus on the cubic NLS equation, which, in its most general form, reads 
\begin{align} \label{cubicNLS}
    i \partial_t + \Delta u - Vu = \kappa \vert u \vert^2 u, \ \ \kappa = \pm 1,
\end{align}
where $V$ denotes a potential. Given the ubiquity of these equations in fields such as nonlinear optics, \cite{ChiaoGarmireTownes}; plasma physics \cite{zakharov1972collapse}; hydrodynamics \cite{Zakharov1968}, obtaining robust numerical simulations is a central problem in mathematical physics.

In this section, we study three different NLS problems of increasing difficulty: 1d NLS with a double well potential, Ginzburg-Landau vortices, and excited states for NLS. By combining the high-precision PINNs framework from Wang {\it et al.}~(2025) \cite{WangEtAl2025Discovery} with the gradient-normalized residual approach, we discover a series of new unstable solutions to each problem, resolved to machine precision. In addition, this approach captures new empirical laws that govern the subtle local features of these solutions. This section will focus on describing the main findings for these three cases. The detailed PINN framework is provided in the supplementary materials.

\begin{figure}[t]
	\centering
	\includegraphics[width=\linewidth]{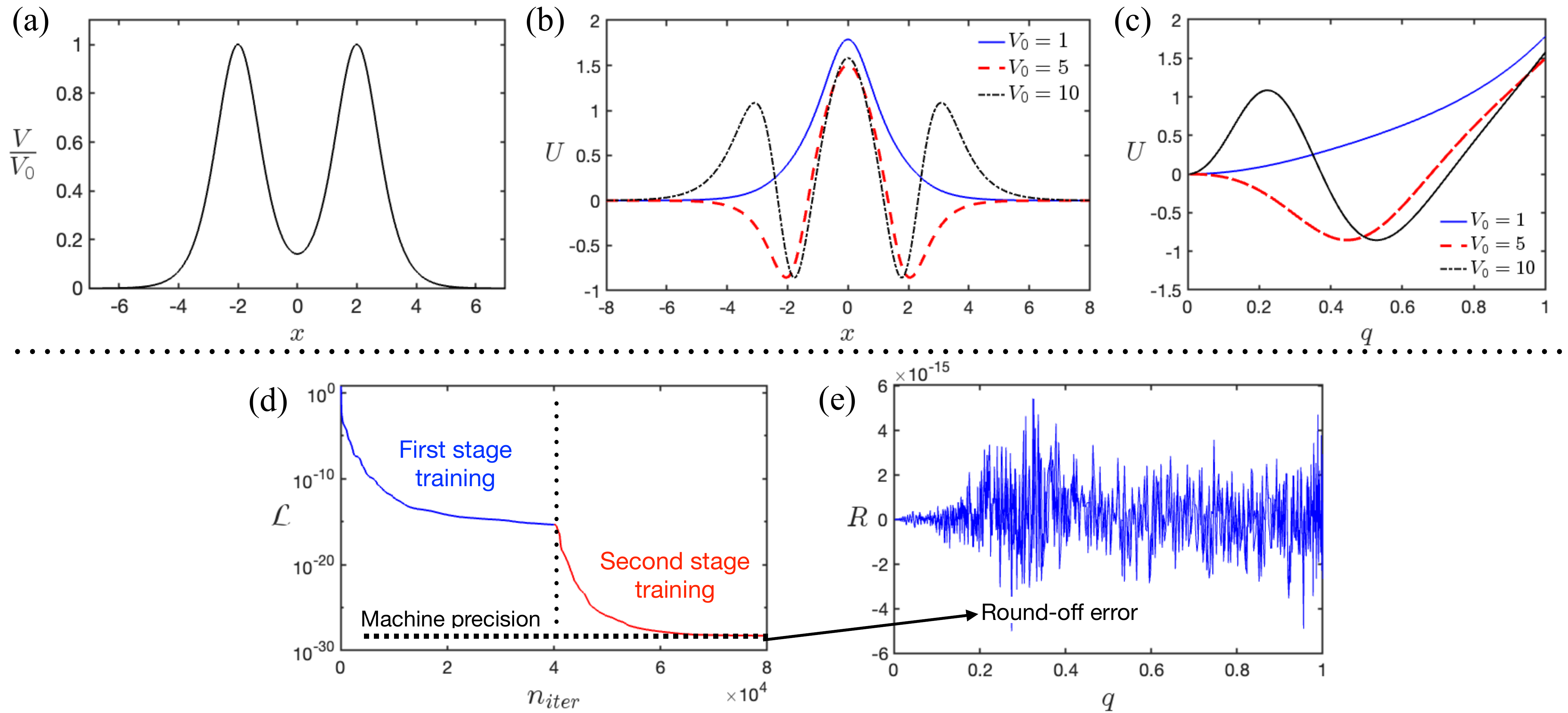}
	\caption{\small \textbf{Solution of the double-well potential problem} (a) The shape of the double-well potential \eqref{eq:doublewell} used in this study width $d=1$ and $x_0=2$. (b) The solution to the Nonlinear Schr\"odinger (NLS) equation \eqref{eq:NLSdouble} with the double-well potential \eqref{eq:doublewell} for different amplitudes $V_0$. (c) The solution profiles in the finite transformed $q$ coordinate, where $q=\sech(\sqrt{\mu}x/2)$. (d) The value of the loss during the training using Gauss-Newton optimizer and multistage scheme. The loss, based on the mean squared error of the equation finally reaches a plateau at $10^{-29}$, which achieve the double-float point machine precision. (e) The PDE residual of the solution for the potential amplitude $V_0=10$. The random spikes show the round-off error of the double-float precision.}
	\label{fig:doublewell}
\end{figure}

\subsection{NLS with Double-well Potential} \label{sec:NLSdoublewell}
The first setting we consider is motivated by nonlinear optics, where Equation \eqref{cubicNLS} in one space dimension and $\kappa = -1$ is a canonical model for transmission in fiber optics \cite{KarpmanKrushkal1969}. The signals being transmitted are both time periodic and spatially localized, which mathematically motivates the ansatz $u(t,x) = \phi(x) e^{i \mu t}, \mu \in \mathbb{R},$ with $\phi$ satisfying $\lim_{x \to \pm \infty} \phi (x) = 0.$ For concreteness, we consider the double-well potential,
\begin{equation}\label{eq:doublewell}
    V(x) = V_0 \left( \sech^2 \left(\frac{x-x_0}{d}\right) + \sech^2 \left(\frac{x+x_0}{d}\right) \right)\, ,
\end{equation}
where $V_0$, $d$, and $x_0$ represent the amplitude, width, and location of the well, respectively. Without loss of generality, in this study, we consistently take $d = 1$ and $x_0 = 2$ (Fig.~\ref{fig:doublewell}a). Then $\phi$ satisfies the boundary value problem
\begin{align}\label{eq:NLSdouble}
\begin{cases}
    \phi'' + \phi^3 - V(x) \phi = \mu \phi \, ,\\ 
    \lim_{x \to \pm \infty} \phi (x) = 0 \, .
\end{cases} .
\end{align}
Moreover, we focus on centered solitons, which are unstable and exhibit even symmetry with respect to the origin in this setting. The key features of the solution vary with the parameter $\mu$; for example, given a potential well amplitude $V_0$, the solution exhibits a double hump for small $\mu$, but a single peak for large $\mu$. Since the problem is also set on an unbounded domain, and that the solution decays exponentially like $e^{-\sqrt{\mu} x}$, the PINN framework from Wang {\it et al.}(2025) \cite{WangEtAl2025Discovery} is well-suited to this situation. We introduce the new coordinate $q = \sech(\sqrt{\mu} x/2)$ (Fig.~\ref{fig:doublewell}c) which maps the original unbounded domain into a finite domain and enforces the even symmetry of the solution. Using the Gauss-Newton optimizer along with the multistage training scheme from Wang \& Lai (2024) \cite{Wang-Lai:multistage-nn}, we successfully recovered these known solutions (Fig.~\ref{fig:doublewell}b), achieving residual errors down to double-float machine precision (Fig.~\ref{fig:doublewell}e) after two stages of training (Fig.~\ref{fig:doublewell}d). We note that this problem does not necessarily require the use of the gradient-normalized residual. Rather, it serves as a test case for applying the PINN method to the NLS equation.

\subsection{Gross-Pitaevskii Vortices}\label{sec:gpvortices}
Equation \eqref{cubicNLS} is also the master equation in the theory of superfluidity \cite{Gross1961, Pitaevskii1961}, where the field $\phi$ is such that $\phi(t,x) e^{it}$ solves \eqref{cubicNLS} with $V=0$ and $\kappa  = 1$, and thus satisfies the Gross-Pitaevskii equation
\begin{align} \label{GL equation}
i \partial_t \phi + \Delta \phi = \phi(\vert \phi \vert^2-1).
\end{align}
It was recognized early on~\cite{Abrikosov} that quantized vortices—special solutions to these equations—could be used to model and explain many experimental observations.

\begin{figure}[t]
	\centering
	\includegraphics[width=\linewidth]{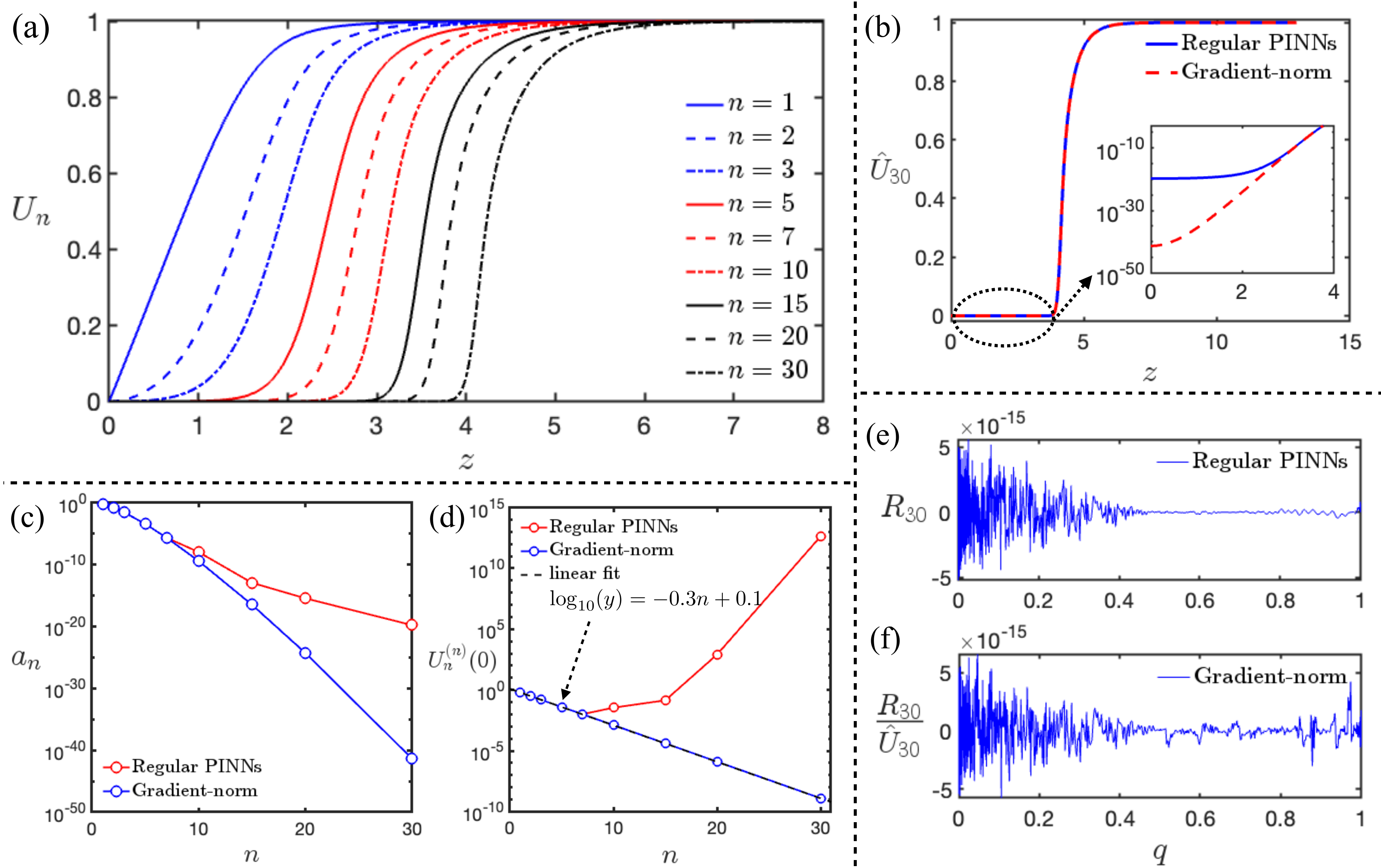}
	\caption{\small \textbf{Solutions of GP-Vortices} (a) spatial profile for GP-vortices solution of different orders from $n=1$ to $n=30$. Here $z = \textrm{arcsinh }{r}$. (b) The renormalized profile of the 30th GP vortex, defined as $\hat{U}_{30} = U_{30}/\eta_{30}$ where $\eta_{30} = r^{30}/(1+r)^{30}$. The red and blue lines show the profile obtained via PINN frameworks using standard equation loss and gradient-normalized loss, respectively. The inset shows both profiles near the origin in a semi-log plot.  (c) The relation of $a_n$ with respect to $n$ extracted from the solutions using two different PINN frameworks. Here $a_n$ refers to the coefficient of the first non-vanish terms of the Taylor expansion of the $n$-th solution at the origin. (d) $U_n^{(n)}(0)$ indicates the first non-vanishing derivative of the $n$-th order solution at the origin, which is equal to $n!a_n$. The derivative plotted with $n$ extracted from the solutions using gradient-normalized PINN align perfectly on a semi-log plot, indicating that the first non-vanishing derivative exponentially decays with the increase of the order $n$. (e) The absolute PDE residual for the 30th GP vortex via the PINN framework using standard equation loss. (f) The relative PDE residual (normalized by the non-vanshing profile $\hat{U}_{30}$) for the 30th GP vortex via the PINN framework using gradient-residual loss.}
	\label{fig:GPvortices}
\end{figure}

Mathematically, vortices constitute a countable family of time-independent solutions of the form $U_n(r) e^{in \theta}$ in two spatial dimensions (here $r$ and $\theta$ denote polar coordinates). These solutions connect two distinct states, represented here by the values $0$ and $1.$ Therefore, the profiles $U_n$ satisfy the boundary value problem
\begin{align} \label{GL:Uprofile}
\begin{cases}
\displaystyle U'' + \frac{1}{r} U' - \frac{n^2}{r^2} U + \big(1 - U^2 \big) U = 0, \vspace{0.2em} \\
\displaystyle U(0)=0, \ \ \lim_{r \to \infty} U(r) = 1
\end{cases} .
\end{align}
Understanding the dynamics of these vortices remains a formidable challenge both theoretically and numerically~\cite{Neu1990, Weinstein-Xin, CarlsonMiller1986}. Obtaining robust numerical approximations is paramount for better simulations and rigorous computer-assisted proofs. The G-P vortices problem \eqref{GL:Uprofile} is defined on an infinite domain. It is straightforward to derive the asymptotic behavior of the profile $U_n$ as $r \to 0$ and $r\to\infty$:
\begin{align}\label{eq:unasymp}
U_n(r) & = a_n r^n \left[1 - \frac{r^2}{4n+4} + O(r^4)\right] \quad \text{as} \  \ r \to 0, \\
U_n(r) & = 1 - \frac{n^2}{2 r^2} + O(r^{-4}) \quad \text{as} \  \ r \to + \infty,
\end{align}
where $a_n>0$ is a unique parameter for each $U_n(r)$ that needs to be determined. It is known \cite{Aguareles2011} that there is a unique monotonic solution to \eqref{GL:Uprofile} such that $U_n(r) \sim a_n r^n.$ The primary numerical challenges in solving the GP vortices problem are handling the infinite domain and accurately determining the vortex core structure $a_n$. Note that the technical difficulty increases with the winding number $n$. To address this issue, we leverage the knowledge of the asymptotic behavior of $U_n$ at infinity \eqref{eq:unasymp} and use a new coordinate transformation $q_n = 1/(1+r^2/n^2)$ to map the domain to a finite interval. Then, applying the same PINN framework from the double-well potential problem (without gradient-normalized residual), we resolved the first 30 vortex profiles with PDE residuals down to $O(10^{-15})$ (Fig.~\ref{fig:GPvortices}e).

\subsubsection{The Vanishing Core of High-order Vortex}
Compared with the double-well potential case, computing vortices for high winding number $n$ presents an additional challenge: the solution profile approaches zero very rapidly as $r\to 0$. More precisely, the vortex profile $U_n$ with large $n$ exhibits a high vanishing order at the origin, namely $U_n\sim r^n,$ see \eqref{eq:unasymp}. To remove this pathological behavior, we introduce a new non-vanishing profile $\widehat{U_n} = U_n/\eta_n(r)$, where $\eta_n=r^n/(1+r)^n$. This new profile $U_{n}$ is non-zero at the origin, positive-definite on the whole domain, and preserves the asymptotics at  infinity. Our numerical results (Fig.~\ref{fig:GPvortices}b) show that the new non-vanishing profile $\hat{U}_n$ for $n=30$ remains extremely close to 0 near the origin. This directly indicates that the unknown parameter $a_n$ in \eqref{eq:unasymp} must decay rapidly to zero as the winding number $n$ increases. This was first noticed in \cite{Aguareles2011}, where an integral expression of $a_n$ in terms of Bessel functions is given.

Figure \ref{fig:GPvortices}(c) shows the relationship of $a_n$ with the winding number $n$. Each $a_n$ value on the red line is obtained using the PINN framework with uniform weighting that achieves $O(10^{-15})$ accuracy (Fig.~\ref{fig:GPvortices}e). It shows that $a_n$ decays exponentially with $n$ when $n<10$, while gradually stalling at a plateau for $n > 10$ (Fig.~\ref{fig:GPvortices}b). However, this plateau is a spurious artifact caused by the insensitivity of the standard loss function to the extremely small value of $U_n$ around the origin.  

As discussed in Section \ref{sec:fluidSingularity}, the standard loss function for PINNs involves the mean squared error of the equation with uniform collocation points and equal weighting. For large winding number $n$ this approach is flawed. Even though it reaches machine precision in the region close to infinity (see the round-off error around $q=0$ in Fig.~\ref{fig:GPvortices}e) with a PDE residual of $O(10^{-15})$, this is insufficient for the region near the origin. There, the true magnitude of the non-vanishing profile $\hat{U}_n$ is even smaller than $O(10^{-15})$ (Fig.~\ref{fig:GPvortices}b). To resolve the core structure of vortices towards machine precision, one must instead minimize the relative error--the ratio of PDE residual to the solution's magnitude--to a tolerance of $O(10^{-15})$.

\subsubsection{Achieving Machine Precision for the Core Structure}
To capture the subtle features of high-order GP vortices, two critical modifications to the PINN framework are necessary. First, to properly capture the vanishing core, we reformulate the learning problem. Instead of learning the vanishing profile $U_n$, we train the network to learn the non-vanishing profile $\hat{U}_n$ as defined above. Since the range of $\widehat{U}_n$ is $(0, 1)$, we further introduce the ansatz $\widehat{U_n} = \exp(-\text{NN}(q_n))$, where NN indicates the output of a fully-connected network, and its input uses the transformed coordinate $q_n = 1/(1+r^2/n^2)$ as defined above. This ansatz effectively enforces the even symmetry of the non-vanishing profile $\hat{U}_n$ and tasks the network with learning its negative logarithm, which has a large stable value around the origin. It also avoids the numerical difficulties of optimizing functions near zero.

Second, we employ the gradient-normalized residual approach with the detailed implementation described in Section 4. Combining these two modifications with the Gauss-Newton optimizer and multistage training allows us to reduce the relative PDE residual $R_{30}/\hat{U}_{30}$ to machine precision $O(10^{-15})$ across the entire domain. Figure \ref{fig:GPvortices}(f) confirms that the relative residual reaches the round-off error in the transformed $q$ coordinates. The resulting non-vanishing profile for the $30$-th vortex is shown to be significantly smaller near the origin than the spurious solution found by the PINN framework using the standard loss function (inset of Fig.~\ref{fig:GPvortices}b).

This unprecedented level of precision allows us to extract the core parameter $a_n$ for each vortex, up to $n=30$. We find that the $a_n$ value decays faster than exponentially. Our results agree with the integral representation given in \cite{Aguareles2011}. We also significantly improve their numerical result, which covered $n$ up to $11.$ 
To uncover the precise numerical relationship between $a_n$ and $n$, we analyzed the first non-vanishing derivative $U_n^{(n)}(0) = n!a_n$ as a function of $n$ and discovered that the data form a precise straight line on a semi-log plot  (Fig.~\ref{fig:GPvortices}d). Based on linear regression, we deduce a novel empirical expression for the core parameter $a_n = \exp(-0.3n+0.1)/n!$. This is inaccessible to the regular PINN approach and, to the best of our knowledge, represents a new contribution to the literature.

\subsection{Excited States for NLS}\label{sec:NLSexcited}

For the final case, we consider the excited states of the 2D cubic NLS equation with $V=0$ and $\kappa = 1$ in \eqref{cubicNLS}. We search for localized solutions of the form $e^{it} U(r) e^{in\theta}$. The radial profile $U(r)$ is governed by the boundary value problem:
\begin{align} \label{NLS:Uprofile}
\begin{cases}
\displaystyle U'' + \frac{1}{r} U' - \frac{n^2}{r^2} U - U+U^3 = 0  \vspace{0.2em}\\
\displaystyle U(0)=0, \quad \lim_{r \to +\infty} U(r) = 0
\end{cases} .
\end{align}
The solution's behavior near the origin is similar to G-P vortices. Thus we will denote $U_n$ the profiles such that $U_n(r) \sim b_n r^n$ as $r\to 0$, where $b_n$ is a unique, unknown parameter to be determined. However, unlike the GP-vortices, the NLS excited states decay to zero at infinity. Its profile is thus dominated by the linear part of the equation \eqref{eq:NLSdouble}, a modified Bessel equation, which drives the asymptotic behavior of the solution profile at infinity:
\begin{equation}\label{eq:excitedDecay}
    U_n(r) \sim r^{-1/2} e^{-r} \quad \text{as} \ r\to\infty.
\end{equation}

\begin{figure}[t]
	\centering
	\includegraphics[width=\linewidth]{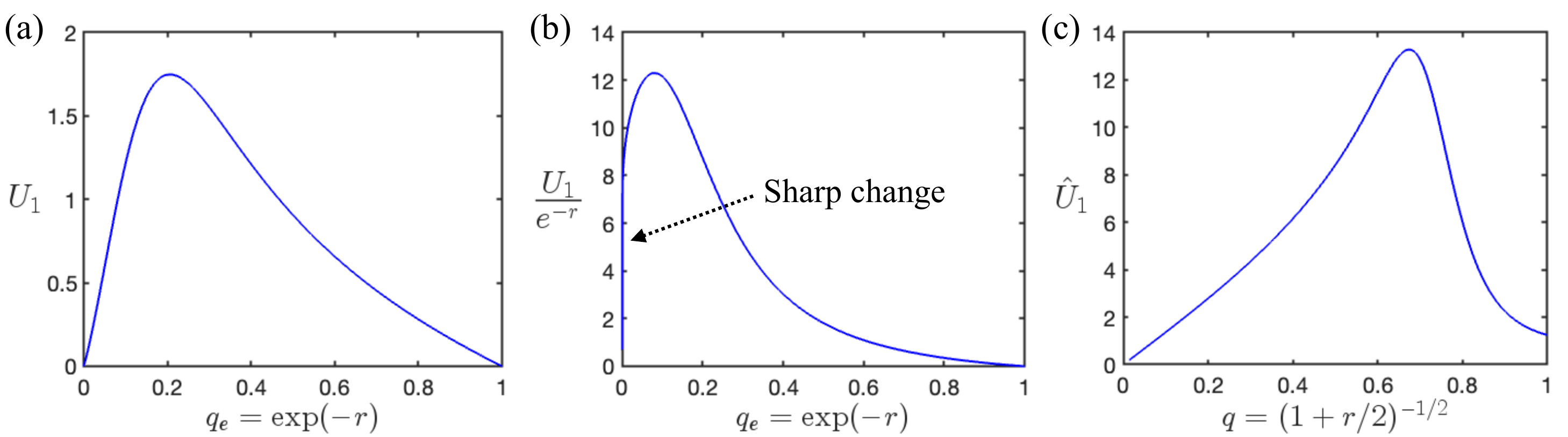}
	\caption{\small \textbf{Coordinate selection of excited solutions} (a) The profile of the 1st NLS excited solution in the transformed coordinate $q_e = \exp(-r)$. (b) The solution normalized by $\exp(-r)$ in the $q_e$ coordinate exhibits a sharp change around $q=0$ (i.e. $r\to\infty$). (c) The non-vanishing profile $\hat{U}_1$ has a regular shape with no sharp change in the new $q=(1+r/2)^{-1/2}$ coordinate. }
	\label{fig:excitedcoord}
\end{figure}

\subsubsection{Subtle Behaviors at Infinity}
First, we focus on solutions $U_n$ to \eqref{NLS:Uprofile} that are positive. In this case, uniqueness is expected. The fast decay at infinity given by \eqref{eq:excitedDecay} is the main issue to address. A basic coordinate transformation tailored to the dominant exponential decay, such as $q_e = \exp(-r)$, is insufficient. Although it may appear to capture the rough shape of the solution (Fig.~\ref{fig:excitedcoord}a), it completely misses the slow power-law decay at infinity. This is illustrated in Fig.~\ref{fig:excitedcoord}(b) where, after normalizing by the expected exponential decay and plotting the solution in the new coordinate $q_e = \exp(-r)$, we observe that the entire asymptotic region where the behavior now follows the power-law decay ($r^{-1/2}$) shrinks into a near-vertical line. This indicates that the subleading power-law decay cannot be captured by this coordinate change. Therefore, a more sophisticated solution ansatz with an improved coordinate transformation is required.

To overcome this challenge, we directly learn the solution profile after normalizing by the dominant exponential decay $\exp(-r)$. Furthermore, since high-order solutions satisfy $U_n \sim b_n r^n$ at the origin, we adopt the same strategy as for the GP-vortices and define a non-vanishing profile 
\begin{equation}\label{eq:excitednonUn}
    \widehat{U}_{n} = \frac{U_n}{\eta_n \cdot e^{-r}} \qquad \text{where} \quad \eta_n = r^n/(1+r)^n,
\end{equation}
which factors out both the high-order vanishing at the origin and the exponential decay at infinity. The resulting profile is now non-zero at the origin and exhibits a simple power-law decay ($r^{-1/2}$) at infinity. Given this simpler asymptotic behavior, we introduce the new coordinate $q = (1+r/2)^{-1/2}$ which maps the infinite domain to the interval $(0, 1]$. 
Figure \ref{fig:excitedcoord}(c) shows the resulting profile in the new coordinate $q.$ It is smooth and well-behaved, demonstrating the success of our combined ansatz and coordinate transformation in capturing the subtle features of the solution at infinity.

\begin{figure}[t]
	\centering
	\includegraphics[width=\linewidth]{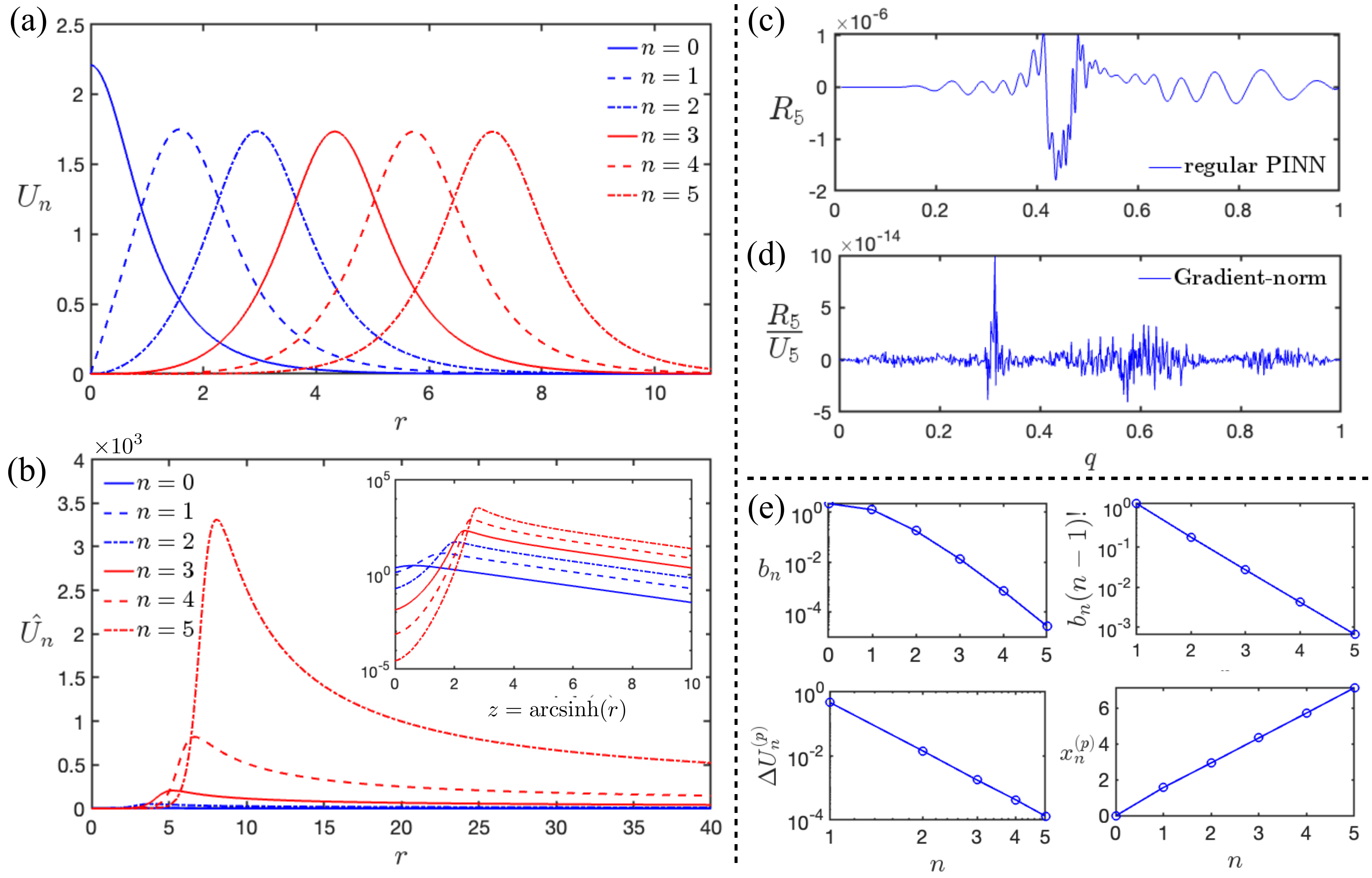}
	\caption{\small \textbf{Solution of NLS excited states} (a) Radial profiles $U_n$ of the NLS excited states for different $n$. (b) Non-zero profile $\hat{U}_n$ \eqref{eq:excitednonUn} of the excited states for different $n$. The inset shows the profiles in a log scale. (c) The absolute PDE residual for the 5th excited states using the PINN framework with standard equation loss, which stagnate at $O(10^{-6})$. (d) The relative PDE residual for the 5th excited states using gradient-normalized residual, combined with Gauss-Newton optimizer and multistage training, which reaches machine precision. (e) Empirical relations for key features of high-order excited states with respect to $n$, including the core parameter $b_n$, the peak magnitude $U_n^{(p)}$, and the peak position $x_n^{(p)}$.}
\label{fig:NLSexcited}
\end{figure}

\subsubsection{Resolving the High-gradient Non-vanishing Profile}
We next focus on positive solutions $U_n$ to \eqref{NLS:Uprofile} for large $n$. While the non-vanishing $\widehat{U}_n$ (which normalizes $U_n$ by the exponential decay) successfully handles the behavior at infinity, it introduces a new numerical challenge: an extremely high gradient.  This pathology arises because of the simple shift of the solution's peak in $U_n(r)$ with $n$ (Fig.~\ref{fig:NLSexcited}a), which results in a massive amplification in $\widehat{U}_n$ (Fig.~\ref{fig:NLSexcited}b). The division by the decay $e^{-r}$—which is vanishingly small at the location of the shifted peak—is the direct cause of this large gradient. In addition, consistent with G-P vortices, the non-vanishing profile $\widehat{U}_n$ of the excited state also has an extremely low value near the origin for high $n$. This creates a massive dynamic range with the peak-to-origin magnitude ratio reaching nearly $O(10^8)$ for $n=5$ (inset of Fig.~\ref{fig:NLSexcited}b).

This large gradient of the non-vanishing profile causes the PINN training with standard equation loss to stagnate (Fig.~\ref{fig:NLSexcited}c), mirroring the behavior seen with the CCF equation. Concurrently, the extremely small value of $\widehat{U}_n(r)$ near the origin makes it difficult to accurately compute the core parameter $b_n$, presenting the same challenge encountered with the GP-vortices. The combination of these two issues significantly increases the difficulty of obtaining high-precision solutions for high-order excited states.

Fortunately, the enhanced PINN framework developed in this work is perfectly suited to overcome both issues. By analogy with GP vortices, we employ two key remedies. First, by introducing the solution ansatz $\widehat{U_n} = \exp(-\text{NN}(q))$, the network learns the negative logarithms of the non-vanishing profile. This significantly reduces the large dynamic range of the solution to a well-behaved target for the network. Second, we employ the gradient-normalized loss function, which ensures that the optimization remains sensitive to the solution's features in both high-gradient and low-magnitude regions. Combining these approaches with a Gauss-Newton optimizer and multistage training, we successfully resolved the first five excited states to machine precision, as confirmed by the relative PDE residual reaching round-off errors throughout the domain (Fig.~\ref{fig:NLSexcited}d). 

This unprecedented level of precision enabled the discovery of several novel empirical relations for the key features of the high-order excited states with respect to $n$. First, the core parameter $b_n$ decays faster than exponentially. However, the empirical relation is slightly different from that of GP vortices. We find that it is $U^{(n)}_n(0)/n = (n-1)!b_n$, rather than $U^{(n)}_n(0)$, that varies linearly with $n$ on a semi-log plot. Second, the solution's peak value, $U^{(p)}_n$, decays monotonically. More precisely, the quantity $\Delta U^{(p)}_n = U^{(p)}_{n-1} - U^{(p)}_{n}$ follows a power law $n$ with a best-fit exponent around 5, namely $\Delta U^{(p)}_n \propto n^{-5}$. Lastly, the position of the peak increases linearly with $n$ with a best-fit slope of approximately 1.37, that is $r_{\text{peak}} \approx 1.37n$. To the best of our knowledge, these high-fidelity relationships are new and provide crucial predictive insights into the behavior of higher-order solutions.

\begin{figure}[t]
	\centering
	\includegraphics[width=\linewidth]{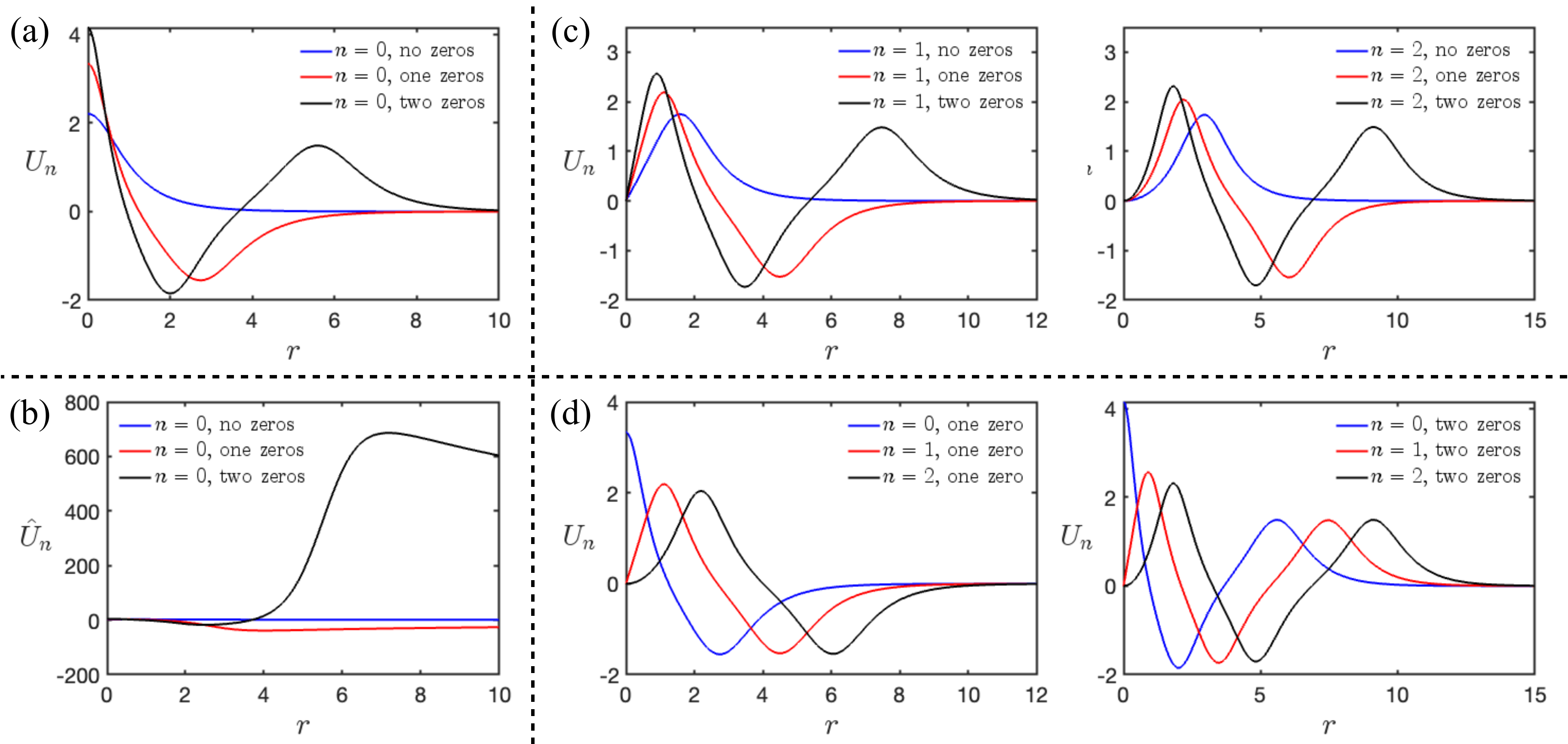}
	\caption{\small \textbf{Excited state solution with more zeros} (a) The radial profiles of NLS excited states for $n=0$ with different nodes (zeros). (b) The non-vanishing profiles $\hat{U}_n(r)$ \eqref{eq:excitednonUn} of the excited states for $n=0$ with diifferent nodes. (c) The radial profile of excited states with different nodes for $n=1$ and $n=2$. (d) Comparison of the radial profiles for a fixed number of nodes with different $n$. }
	\label{fig:excitedzeros}
\end{figure}

\subsubsection{Non-uniqueness of the Excited State Solution}
We now come to profiles $U_n$ that change sign. Indeed, an additional difficulty in the NLS excited states problem \eqref{NLS:Uprofile} is caused by non-uniqueness. The vanishing order $n$ of $U_n$ at $r=0$ does not uniquely characterize the solution. There exists a countable family of solutions, indexed by their number of nodes (zeros) $k$, for any given $n$. Uniqueness is only expected for a given pair $(n,k)$ \cite{Berestycki-Lions:nonlinear-scalar-field-i-ground-state,Berestycki-Lions:nonlinear-scalar-field-ii-existence-infinitely-solutions}. These nodal solutions ($k>0$) are inherently unstable, and to the best of our knowledge, none of them have been numerically discovered in the literature prior to this work. The primary challenge in finding these states lies in locating the nodes. Without additional guidance, the training process naturally converges to the nodeless $(k=0)$ solution, as neural networks exhibit a spectral bias towards the lowest-variation solution \cite{xu2019frequency, rahaman2019spectral}.

To overcome this issue, we introduce targeted inequality constraints to guide the training into the desired solution basin. For example, to find the $(n, k)=(0, 1)$ state, we leverage the fact that its profile must cross the horizontal axis and then decay to zero while remaining negative. By imposing the inequality constraint $U_0(r_c) < -0.1$ at a chosen control point $r_c$ (e.g., $r_c=5$), we successfully identify the desired single-node solution for $n=0$ (Fig.~\ref{fig:excitedzeros}a). This strategy can be extended to solutions with multiple zeros.

Using this approach, we further discovered the solutions with one and two nodes for $n=1$ and $n=2$ (Fig.~\ref{fig:excitedzeros}c). For $k=2,$, it was required to impose multiple inequality constraints. Comparing the radial profile for a fixed number of nodes $(k>0)$, we find that both the node and peak positions increase linearly with $n$ (Fig.~\ref{fig:excitedzeros}d), consistent with the pattern observed for the $k=0$ solutions (Fig.\ref{fig:NLSexcited}e). This provides a powerful predictive tool for estimating the structure of high-order nodal solutions. We note that the solution becomes increasingly oscillatory as $k$ gets larger, which also induces high gradients in the corresponding non-vanishing profiles $\hat{U_n}$ \eqref{eq:excitednonUn} (Fig.~\ref{fig:excitedzeros}b). Thus, the enhanced PINN framework, using gradient-normalized residuals, remains essential to resolve these solutions with high precision. However, a limitation arises: since the profiles are no longer positive-definite, the logarithmic ansatz $\hat{U}_n = \exp(-\text{NN}(q))$ is no longer applicable. Although a standard network can resolve the profiles for low $n$ and $k$, the design of a more sophisticated ansatz to handle the oscillations, subtle core structure around the origin, and high-gradient features at infinity for solutions with large $n$ and $k$ remains a compelling direction for future research.

\section{Method}

The present paper builds upon the computational frameworks for discovering unstable singularities established in~\cite{Wang-Lai-GomezSerrano-Buckmaster:pinn-selfsimilar-boussinesq,WangEtAl2025Discovery}. These frameworks combine mathematical analysis with advanced Physics-Informed Neural Networks (PINNs) to solve various nonlinear partial differential equations. Although this approach has robustly identified unstable solutions for diverse fluid equations, the accuracy achieved is highly dependent on the solution's stability, often struggling with the large gradients inherent in higher-order unstable modes. In this section, we recap the key techniques used in the foundational framework and discuss the detailed implementation of the gradient-normalized residue into the training process.

\subsection{Foundational frameworks of singularity discovery}\label{sec:basicFW}

The core principle of the framework developed in~\cite{Wang-Lai-GomezSerrano-Buckmaster:pinn-selfsimilar-boussinesq,WangEtAl2025Discovery} is the embedding of mathematical insights directly into the neural network architecture. This transforms a generic function approximation task into a highly constrained search for mathematically relevant solutions. This framework consists of the following key components:

{\bf Self-similar coordinates}:  Rather than solving the time-dependent fluid equations, we introduce self-similar coordinates around the singularity and derive the corresponding steady-state self-similar equation, which is defined over an infinite spatial domain in the local coordinates.

{\bf Embedding mathematical insights}: To solve the transformed equation, we leverage the asymptotic behavior of the solution to construct a coordinate transformation that maps the infinite domain to a normalized finite domain. Furthermore, we enforce known mathematical properties, such as symmetries, asymptotic decay rates, and boundary conditions, by hard-coding them into a structured solution \emph{ansatz}. This ansatz shapes the network's output, ensuring that the fundamental constraints are satisfied by construction.

{\bf Multi-stage neural network}: To enhance the expressiveness of the neural network and reduce the PDE residuals, we implement the Multi-Stage Neural Network (MSNN) \cite{Wang-Lai:multistage-nn}, where the residual of the first neural network is used to optimize the design of the subsequent neural networks.

{\bf Second-order optimizer}: To train the PINNs, we eschew standard first-order optimizers in favor of a powerful second-order optimizer utilizing a full-matrix Gauss-Newton method \cite{martens2015optimizing}. Implemented via \texttt{kfac-jax}~\cite{kfac-jax2022github}, this approach includes an automated scheme to determine the learning rate and momentum coefficient. This optimizer yields significantly faster convergence than standard first-order methods without requiring extensive parameter sweeps. Other forms of Gauss-Newton optimization has also been previously implemented in PINNs, e.g., \cite{pmlr-v202-muller23b,jnini2024gaussnewtonnaturalgradientdescent, guzmáncordero2025improvingenergynaturalgradient}.

{\bf Collocation point sampling}: We employ a hybrid sampling strategy to select collocation points. The first method is \emph{location-based sampling}, where points are drawn from fixed spatial regions following a user-prescribed distribution. The second is \emph{adaptive sampling}, e.g., \cite{doi:10.1137/19M1274067,https://doi.org/10.1111/mice.12685, WU2023115671}). In our implementation collocation points are sampled globally but weighted by the square of the PDE residuals. This concentrates the training effort on regions where the current PDE residuals are highest.

{\bf $\lambda$ identification}: The final component is the identification of the admissible self-similar scaling parameter $\lambda$. This is accomplished through two complementary strategies: (1) direct \emph{analytical inference}, where $\lambda$ is derived from regularity conditions at the origin , and (2) an iterative \emph{funnel inference} method that utilizes a secant search to find the $\lambda$ that minimizes PDE residuals. Regardless of the method employed, the objective is to identify the $\lambda$ that minimizes the maximum PDE residual both near the origin and across the entire domain.

\subsection{Resolving the sharp gradient by gradient normalization}\label{sec:grad_method}

The primary obstacle for the foundational framework described above to discover new, more unstable solutions is the emergence of extreme, localized gradients in the solution profiles. As the unstable mode increases, these gradients intensify rapidly (Fig.~\ref{fig:high-gradient}), creating a pathological optimization landscape that fundamentally undermines PINN training strategies.

\subsubsection{Limitation of standard $L^2$ loss formulations}
In standard PINN formulations, the loss function typically comprises the mean squared error of both the solution constraints (i.e., boundary and initial conditions) and the PDE residuals. However, utilizing the framework described above, the solution constraints are strictly enforced via the hard-coded ansatz. Consequently, the loss function reduces to the weighted $L^2$ norm of the PDE residuals evaluated at collocation points across the computational domain:
\begin{equation}\label{eqn:loss}
\mathcal{L} = \frac{1}{N_d}\sum_{k=1}^{N_e}\sum_{i=1}^{N_d} w_i |\mathcal{R}_k(y_i)|^2,
\end{equation}
where $\mathcal{R}_k(y_i)$ denotes the residual of the $k$-th governing equation evaluated at the collocation point $y_i$. $N_e$ represents the total number of coupled equations, $N_d$ is the number of collocation points in the batch, and $w_i$ represents the weight assigned to the residual at point $y_i$. 

However, for partial differential equations involving solutions with sharp gradients, the magnitude of the PDE residual inherently scales with the local gradient of the solution. Consequently, the numerical cancellation required to satisfy $\mathcal{R}(y) \approx 0$ becomes increasingly sensitive to small approximation errors. In this scenario, the adaptive sampling strategy used in the foundational framework responds by concentrating collocation points (or equivalently increasing weights $w_i$) in these high-gradient regions. While intuitively sound, we observe that this approach causes the training to plateau, preventing the PDE residual from being minimized further (Fig.\ref{fig:ccf_issue}b and \ref{fig:IPM4th}b). 

\subsubsection{Theoretical Motivation: Mimicking an Adaptive Coordinate System}

The fundamental challenge posed by high-gradient solutions is that standard loss functions prioritize minimizing absolute errors. This disproportionately weights high-gradient regions, where residuals are inherently larger, often leading to optimization stagnation  while obscuring subtle features elsewhere in the domain.

Conceptually, this pathology could be resolved by an idealized auxiliary coordinate transformation that stretches high-gradient regions and compresses low-gradient regions, such that the solution exhibits uniform variation in the new coordinates. In such an idealized system, the optimization landscape would be well-behaved.

However, explicitly constructing such a transformation is often intractable that could lead to a severely distorted coordinate system, potentially rendering the transformed PDE more difficult to solve.

Instead of explicitly changing coordinates, we aim to mimic the regularization effects of this idealized transformation by modifying the sampling strategy and the loss function directly in the original coordinates. This requires two key adjustments.

First, uniform sampling in the idealized coordinates corresponds to non-uniform sampling in the physical coordinates, with density proportional to the Jacobian of the idealized transformation. Since the goal is to uniformize the gradient, this idealized Jacobian must scale proportionally to the solution's gradient magnitude. This provides a theoretical justification for concentrating collocation points in high-gradient regions.

Second, and crucially, we must rebalance the optimization. To ensure the optimization behaves as if it were operating in the idealized coordinates, the physical residual must be renormalized by the same scaling factor -- the proxy for the idealized Jacobian.

This analysis reveals that merely increasing sampling density is insufficient. To robustly solve high-gradient problems, one must simultaneously normalize the PDE residual by the local gradient magnitude. This strategy rebalances the optimization landscape, ensuring the training process targets uniform \emph{relative} error rather than uniform \emph{absolute} error.

\subsubsection{Implementation of gradient normalization}
Recall that in our foundational framework, we employ coordinate compactification to map infinite domains onto normalized finite domains. For the problems studied here, the dependent variables for the IPM equation (specifically the density $H$) and the Gross-Pitaevskii (G-P) vortices are monotonic functions in these compactified coordinates. Consequently, we can directly leverage the derivative of the solution to normalize the governing equations, following the theoretical derivation in the previous section.

\paragraph{Monotonic Case: IPM and G-P vortices.}
We take the IPM problem (Section \ref{sec:ipm}) as a primary example. The critical quantity governed by the nonlinearity is the density field $H$. In the compactified coordinate $q$, $H$ is monotonic. Its spatial gradient corresponds physically to the vorticity, $\Omega$, which is non-negative (vanishing only at the origin and infinity). 

However, direct normalization by $\Omega$ is numerically unstable due to its zeros at the boundaries. To circumvent this, we treat $\Omega$ as an auxiliary output constrained by the relation $\partial_{y_1}H = -\Omega$ and enforce a structured ansatz:
\begin{equation}\label{eq:vorticity}
    \Omega(y) = -\frac{y_1}{\sqrt{1+y_1^2+y_2^2}}\cdot \exp(\text{NN}_\Omega(y)) \cdot q(y_1, y_2),
\end{equation}
where $q(y_1, y_2)$ captures the asymptotic decay at infinity. This construction ensures that $\Omega$ maintains the correct sign and boundary behaviors. Crucially, the term $\exp(\text{NN}_\Omega)$ represents the core gradient magnitude and is \textit{strictly positive} everywhere. Thus, instead of normalizing by the full vorticity $\Omega$, we normalize the PDE residual by $\exp(\text{NN}_\Omega)$.

A practical challenge is that $\Omega$ (and thus $\text{NN}_\Omega$) evolves during training. We explored two strategies to address this: (1) \textit{Pre-training}, where a preliminary run identifies the rough profile of $\exp(\text{NN}_\Omega)$ to be used as a fixed normalization factor; and (2) \textit{Dynamic updates}, where the normalization factor is updated periodically (e.g., every 500 iterations) based on the current network state. Systematic testing indicates that both methods yield comparable accuracy. Adopting the latter, the modified loss function is:
\begin{equation}
\mathcal{L}_{\text{norm}}^{\text{(IPM)}} = \frac{1}{N_d}\sum_{k=1}^{N_e}\sum_{i=1}^{N_d} \left| \frac{\mathcal{R}_k(y_i)}{\exp(\text{NN}_\Omega(y_i))}\right|^2.
\end{equation}
Since $\exp(\text{NN}_\Omega)$ captures the solution's gradient feature, we utilize it directly as the probability density function (PDF) for sampling collocation points, naturally concentrating computational effort in high-gradient regions. This approach is also applied to the G-P vortices problem (Section \ref{sec:gpvortices}), which exhibits similar monotonicity (Fig.~\ref{fig:GPvortices}a).

\paragraph{Non-Monotonic Case: CCF and NLS.}
Unlike IPM, the CCF equation (Section \ref{sec:ccf}) and the 2-D cubic NLS equation (Section \ref{sec:NLSexcited}) govern variables (such as vorticity $\Omega$ in CCF and excited states $U_n$ in NLS) that are non-monotonic (Fig.~\ref{fig:high-gradient}a and Fig.~\ref{fig:NLSexcited}a). 

However, the solution ansatz for the relevant physical variables in these problems is similar to \eqref{eq:vorticity}, and employs an exponential core to capture multiscale variations. Therefore, we can generalize the normalization strategy by approximating the gradient magnitude as $\exp(\alpha \cdot \text{NN}_\Omega)$, where $\alpha$ is a tunable parameter balancing the normalization factor with the scale of the solution's gradient. For the CCF equation, we find $\alpha \approx 2$ to be effective. The resulting loss function is:
\begin{equation}
\mathcal{L}_{\text{norm}}^{\text{(CCF)}} = \frac{1}{N_d}\sum_{k=1}^{N_e}\sum_{i=1}^{N_d} \left| \frac{\mathcal{R}_k(y_i)}{\exp(2 \cdot \text{NN}_\Omega(y_i))}\right|^2.
\end{equation}
Furthermore, when employing Gradient-Enhanced PINNs (g-PINNs) \cite{yu2022gradient} -- which include the derivative of the PDE residuals in the loss—the residuals involve higher-order derivatives of the solution and thus exhibit even steeper gradients. To properly normalize these terms, we apply the same exponential factor but with a larger scaling parameter $\alpha$.

Specific implementation details, including network architectures and hyperparameter choices for each problem, are provided in the Supplementary Information.

\section{Conclusion}
Computing highly unstable solutions to nonlinear PDEs presents a fundamental optimization pathology: standard $L^2$ residual minimization creates a loss landscape dominated by high-gradient regions, effectively obscuring the subtle signals that determine critical physical parameters such as the self-similar scaling rate $\lambda$. We have resolved this pathology through adaptive gradient-based residual normalization, which reconditions the loss function to measure solution quality relative to local gradient magnitude rather than the absolute value of the residual. This rebalancing, combined with multi-stage correction networks employing specialized high-frequency architectures, enables systematic convergence to machine precision, even for extremely unstable phenomena characterized by exponentially sharp gradients.

This framework achieves residuals below $O(10^{-13})$ for all known unstable self-similar solutions in the CCF equation and between $O(10^{-11})$ and $O(10^{-13})$ for the 2-D IPM equations, including confirmation of the fourth unstable IPM solution that was previously inaccessible due to insufficient accuracy. With the high-precision solution achieved, we validate the empirical relation between the non-smooth signal at the origin and the accuracy of the scaling parameter $\lambda$ for different orders of unstable solutions, and confirm that a high-precision framework capable of resolving sharp local gradients is indispensable for discovering highly unstable singularities.

Moreover, we demonstrate that this framework can be successfully applied to resolve other challenging problems in mathematical physics. Applied to NLS equations, the methodology achieves machine precision for unstable solitons, Gross-Pitaevskii vortices with winding numbers up to $n=30$, as well as new excited states for the NLS equation. Furthermore, the method unveils subtle local features that were previously unknown in the literature. For instance, we identified an empirical law for the vortex core coefficient, $a_n \propto e^{-cn}$, and observed that the peak locations of NLS excited states scale linearly with the winding number.

Beyond numerical discovery, high accuracy is a prerequisite for rigorous mathematical verification via computer-assisted proofs.  The adaptive normalization strategy is general: it applies in settings where gradient pathologies obstruct optimization. By embedding mathematical structure into the learning process to handle ill-conditioned numerics, this work offers a robust methodology for problems where both traditional numerical methods and generic neural network approaches are ineffective.

\section{Acknowledgments}

C-Y. Lai acknowledges the NSF  via grant  DMS-2245228. T. Buckmaster, T. L\'{e}ger and Y. Wang were supported in part by the NSF grant DMS-2243205 and the Simons Foundation Mathematical and Physical Sciences collaborative grant `Wave Turbulence'. T. Buckmaster was also supported by the NSF grant DMS-2244879. Y. Wang's research was additionally supported an NYU Postdoctoral Research and Professional Development Support Grant. Both T. Buckmaster and Y. Wang  gratefully acknowledge Stanford Research Computing and New York University High Performance Computing
for providing computational resources and support.

\bibliographystyle{unsrtnat}
\bibliography{references} 

\end{document}